\documentclass[11pt]{amsart}

\usepackage[T1]{fontenc}
\usepackage{lmodern}
\usepackage{microtype}
\usepackage{amsmath,amssymb,amsthm,mathtools}
\usepackage{booktabs,longtable,array}
\usepackage{enumitem}
\usepackage{hyperref}
\hypersetup{hidelinks,pdftitle={Hyperbolic Polynomial Certificates: Universality, Cubic Exactness, and a Weak--Strong Separation},pdfauthor={Jie Wang}}
\usepackage[nameinlink,noabbrev]{cleveref}
\usepackage{url}
\usepackage{listings}
\newif\ifprintcomputercode
\printcomputercodefalse
\usepackage{geometry}
\setlist{leftmargin=2em,itemsep=2pt,topsep=4pt}
\allowdisplaybreaks

\newtheorem{theorem}{Theorem}[section]
\newtheorem{proposition}[theorem]{Proposition}
\newtheorem{lemma}[theorem]{Lemma}
\newtheorem{corollary}[theorem]{Corollary}
\theoremstyle{definition}
\newtheorem{definition}[theorem]{Definition}
\newtheorem{example}[theorem]{Example}
\theoremstyle{remark}
\newtheorem{remark}[theorem]{Remark}

\newcommand{\R}{\mathbb{R}}
\newcommand{\C}{\mathbb{C}}
\newcommand{\Q}{\mathbb{Q}}
\newcommand{\Sph}{\mathbb{S}}

\newcommand{\Ext}{\operatorname{Ext}}
\newcommand{\Rea}{\operatorname{Re}}
\newcommand{\diag}{\operatorname{diag}}
\newcommand{\tr}{\operatorname{tr}}

\newcommand{\calP}{\mathcal{P}}
\newcommand{\calZ}{\mathcal{Z}}
\newcommand{\calB}{\mathcal{B}}
\newcommand{\calH}{\mathcal{H}}
\newcommand{\calE}{\mathcal{E}}
\newcommand{\calC}{\mathcal{C}}
\newcommand{\calT}{\mathcal{T}}
\newcommand{\Lam}{\Lambda_+}

\newcommand{\eps}{\varepsilon}

\title[Hyperbolic polynomial certificates]{Hyperbolic Polynomial Certificates:\\Universality, Cubic Exactness, and a Weak--Strong Separation}
\author{Qingpo Cao}
\address{State Key Laboratory of Mathematical Sciences, Academy of Mathematics and Systems Science, Chinese Academy of Sciences, Beijing, China}
\email{caoqingpo@amss.ac.cn}
\author{Jie Wang}
\address{State Key Laboratory of Mathematical Sciences, Academy of Mathematics and Systems Science, Chinese Academy of Sciences, Beijing, China}
\email{wangjie212@amss.ac.cn}
\date{\today}

\begin{document}

\begin{abstract}
We study certificates of polynomial nonnegativity based on hyperbolic polynomials and their parameterized B\'ezoutians and establish three results. First, in the unrestricted hyperzout model of Ng and Saunderson, every nonnegative real form of degree $2r$ in $m$ variables admits a hyperzout representation using an auxiliary hyperbolic form of degree $2r$ in $m+1$ variables with parameters $k=1$ and $\mu=r+1$. Second, every homogeneous hyperbolic cubic in five real variables is SOS-hyperbolic in Saunderson's fixed-direction sense: for each hyperbolicity direction $e$, $B_{p,e,v}$ is a matrix SOS for every $v\in\Lam(p,e)$. Third, a general localization proposition turns an interior-SOS central Wronskian and a non-matrix-SOS central B\'ezoutian into a two-direction weak--strong separation. An icosahedral seed gives a rational example of degree $46$ in $43$ variables for which every $W_{u,v}$ on the cone is an SOS, but some $B_{p,u,v}$ is not a matrix SOS. The obstruction persists at nearby interior directions, giving a rational hyperbolicity direction at which the polynomial is weakly SOS-hyperbolic but not SOS-hyperbolic in Saunderson's fixed-direction sense. The construction is supported by exact rational Gram and trace certificates.
\end{abstract}


\maketitle

\section{Introduction}

Certifying global polynomial nonnegativity is a basic problem in real algebraic geometry and polynomial optimization. Sums of squares (SOS) provide a sufficient condition with a semidefinite formulation: after choosing a monomial basis, an SOS identity is equivalent to the existence of a positive semidefinite Gram matrix. This connection underlies the polynomial optimization methods of Lasserre and Parrilo \cite{Lasserre2001,Parrilo2003}. SOS representations do not, however, exhaust the nonnegative polynomials. For forms of positive even degree, Hilbert's classification shows that every nonnegative form is an SOS only in the cases of quadratic forms, forms in at most two variables, and ternary quartics \cite{Hilbert1888}. Blekherman quantified the discrepancy between the two cones: at fixed degree greater than two, their normalized volume comparison becomes increasingly unfavorable to SOS as the number of variables grows \cite{Blekherman2006}. These results motivate certificates that exploit structures beyond a polynomial SOS identity.

Hyperbolic polynomials supply one such structure. They originated in G\aa rding's study of linear partial differential equations and connect real-rootedness with convexity. If a homogeneous polynomial $p$ is hyperbolic with respect to a direction $e$, then the roots of $t\mapsto p(x+te)$ define a spectral calculus analogous to the eigenvalue calculus for symmetric matrices, and the associated hyperbolicity cone is convex \cite{Garding1959}. G\"uler proved that $-\log p$ is a logarithmically homogeneous self-concordant barrier on the cone interior \cite{Guler1997}. Renegar developed hyperbolic programming and derivative relaxations as extensions of semidefinite programming \cite{Renegar2006}. Hyperbolicity thus provides both a source of globally nonnegative polynomials and a conic optimization framework in which to search for their certificates.

The basic spectrahedral examples arise from definite determinantal representations
\[
p(x)=\det\!\left(\sum_{i=1}^n x_iA_i\right),
\qquad
\sum_{i=1}^n e_iA_i\succ0,
\]
whose hyperbolicity cones are the spectrahedra $\{x:\sum_i x_iA_i\succeq0\}$. The classical Lax conjecture for three-variable hyperbolic polynomials was established through the work of Lewis--Parrilo--Ramana \cite{LewisParriloRamana2005} and Helton--Vinnikov \cite{HeltonVinnikov2007}. The generalized Lax question asks whether every hyperbolicity cone admits such a linear matrix inequality description. Positive cases include elementary-symmetric and derivative-relaxation cones \cite{Branden2014,SaundersonParrilo2015}, symmetric hyperbolic cubics \cite{BlekhermanLindbergShu2025}, and all hyperbolic cubics in five variables through Netzer's Clifford construction \cite{Netzer2026}. For smooth hyperbolicity cones, Netzer and Sanyal proved the weaker representation as a spectrahedral shadow, that is, a linear projection of a spectrahedron \cite{NetzerSanyal2015}.

It is essential to distinguish the geometry of a cone from the algebraic properties of a polynomial defining it. Br\"and\'en exhibited hyperbolic polynomials for which no positive power has a definite determinantal representation \cite{Branden2011}. Kummer showed that the hyperbolicity cone of the specialized V\'amos quartic is nevertheless spectrahedral \cite{Kummer2016Vamos}. The same quartic, up to a permutation of variables, fails SOS-hyperbolicity, as shown by Saunderson \cite[Proposition~4.9]{Saunderson2019}. Thus a spectrahedral description of a cone does not, by itself, furnish the matrix SOS certificates associated with its particular defining polynomial. This distinction is central to the cubic result below: the proof uses specific algebraic information in the cone representation.

Parameterized Hermite matrices and B\'ezoutians connect these geometric and algebraic questions. Kummer, Plaumann, and Vinzant expressed cones of interlacers as linear sections of cones of nonnegative polynomials and used SOS relaxations to obtain semidefinite approximations of hyperbolicity cones \cite{KummerPlaumannVinzant2015}. Netzer, Plaumann, and Thom related definite determinantal representations to SOS decompositions of Hermite matrices \cite{NetzerPlaumannThom2013}. Kummer subsequently used B\'ezoutian SOS decompositions to prove that every smooth hyperbolic polynomial has a hyperbolic multiplier whose product with it admits a definite determinantal representation \cite{Kummer2017Bezoutians}. This theorem does not assert the cone-containment condition needed to recover the original cone as the positive semidefinite locus. Clifford-algebraic approaches to these representation problems were developed by Netzer and Thom \cite{NetzerThom2014}.

Saunderson developed a complementary use of parameterized B\'ezoutians: their polynomial compressions certify scalar nonnegativity and lead to searches over hyperbolicity cones \cite{Saunderson2019}. Ng and Saunderson organized the resulting Wronskian and B\'ezoutian certificates into the hyperwron and hyperzout classes and established non-representability results, with an additional degree restriction in the hyperzout setting \cite{NgSaunderson2026}. Here the distinction between a fixed auxiliary hyperbolic polynomial and a variable one matters. With the auxiliary data fixed, one studies the certificates generated by a particular cone; allowing those data to depend on the target form changes the representation problem itself.

These developments lead to three questions. How much expressive power does the unrestricted hyperzout model have? When does the B\'ezoutian of a fixed hyperbolic polynomial admit a matrix SOS decomposition? Can scalar Wronskian SOS conditions hold throughout a hyperbolicity cone while the corresponding matrix conditions fail? Our results address these questions through universality, exactness for five-variable cubics, and a weak--strong separation.

We begin with the unrestricted hyperzout model. Let $\calP_{m,2r}$ denote the cone of nonnegative real forms of degree $2r$ in $m$ variables, and let $\calZ_{m,2r}$ denote the set of forms admitting a hyperzout representation in the sense recalled in \cref{sec:prelim-hyperzout}. The dimension arguments of Ng and Saunderson exclude degree-restricted hyperzout representations for suitable degrees and numbers of variables \cite{NgSaunderson2026}. They leave open whether that obstruction reflects a limitation of the unrestricted model or of its degree regime. The following theorem identifies the role of this distinction.

\begin{theorem}[Universality of unrestricted hyperzout certificates]\label{thm:intro-univ}
For all integers $m,r\ge 1$,
\[
\calZ_{m,2r}=\calP_{m,2r}.
\]
More precisely, every nonzero $f\in\calP_{m,2r}$ has a representation with
\[
n=m+1,\qquad d=2r,\qquad k=1,\qquad \mu=r+1,
\]
with hyperbolicity direction and both B\'ezoutian directions equal to $e=(1,0,\ldots,0)$ and with substitution $\phi(x)=(0,x)$.
\end{theorem}

The construction starts from the nonnegativity of $f$. We choose a homogeneous real-rooted polynomial $h(t,x)$ of degree $r$ whose normalized one-variable fibers have simple real roots. For sufficiently small $\eps>0$,
\[
p(t,x)=h(t,x)^2-\eps f(x)
\]
is hyperbolic. A root-moment factorization of the B\'ezoutian $B(p,p')$ then recovers $f$ by compression with an explicit polynomial vector. A triangularity argument for $B(p,1)$ gives the required degree profile. The construction lies on the degree boundary $dk=2r$; for $r\ge2$ it does not satisfy the stricter inequality $d\le 2\mu-3$ used in the degree-restricted subclass. Thus \cref{thm:intro-univ} complements the obstruction theory of \cite{NgSaunderson2026}: degree-restricted hyperzout certificates do not cover all nonnegative forms in the regimes identified there, whereas the unrestricted class is universal.

The significance of this universality statement is representational. The auxiliary polynomial is allowed to depend on the target form, and its hyperbolicity is established using the assumed nonnegativity of that form. The theorem therefore does not give an efficient nonnegativity test or a universal fixed auxiliary cone. Instead, it shows that degree-restricted non-representability cannot be extended to the unrestricted class without further hypotheses. The explicit degree profile also separates the existence of a certificate from the question of how small or computationally accessible it can be.

For the second question, we first recall the known results. Linear forms have constant nonnegative B\'ezoutians, and every hyperbolic quadratic is SOS-hyperbolic, irrespective of the number of variables \cite[Proposition~4.8]{Saunderson2019}. More generally, a determinantal representation of a positive power of $p$ that is definite at the chosen hyperbolicity direction suffices in any degree \cite[Proposition~4.7]{Saunderson2019}. This criterion, together with the Helton--Vinnikov theorem \cite{HeltonVinnikov2007}, implies that all hyperbolic forms in at most three variables are SOS-hyperbolic in any degree \cite[Proposition~4.8]{Saunderson2019}.

In contrast, the specialized V\'amos quartic mentioned above already gives a four-variable counterexample \cite[Proposition~4.9]{Saunderson2019}. Saunderson's arguments for increasing the degree and number of variables propagate this failure: for every $d\ge4$ and $n\ge4$, there is a degree-$d$ hyperbolic form in $n$ variables that is not SOS-hyperbolic \cite[Propositions~4.10 and~4.13, Theorem~4.14]{Saunderson2019}. Cubics therefore require a separate analysis of how SOS-hyperbolicity depends on dimension.

Within degree three, dimension is decisive. Saunderson's graph construction gives a hyperbolic cubic in $43$ variables that is not SOS-hyperbolic \cite{Saunderson2019}. Blekherman, Kummer, Sanyal, Shu, and Sun later constructed a six-variable hyperbolic cubic with a non-SOS mixed Rayleigh difference \cite[Section~6]{blekherman2023linear}. The five-variable case remained the unresolved cubic case in this comparison \cite{NgSaunderson2026}. Our second result positively settles that case and, in particular, rules out the five-variable case of Conjecture~5.9 in \cite{Saunderson2019}.

\begin{theorem}[SOS-hyperbolicity of five-variable cubics]\label{thm:intro-cubic}
Let $p\in\R[z_0,\ldots,z_4]$ be homogeneous, cubic, and hyperbolic with respect to $e\in\R^5$. Then $p$ is SOS-hyperbolic with respect to $e$ in Saunderson's fixed-direction sense: $B_{p,e,v}(x)$ is a matrix SOS for every $v\in\Lam(p,e)$.
\end{theorem}

% After Netzer's normal-form reduction, the essential case is
% \[
% p_q(t,x)=t^3-3t\|x\|^2+2q(x),\qquad x\in\R^4,
% \]
% where $\|q\|_{\infty}\le1$ on the unit sphere. The canonical symmetric pencil $T_q(x)=\frac16\nabla^2q(x)$ contains more information than the scalar identity defining Clifford realizability and is already present in Netzer's construction. The extracted lemma shows that every $q$ in the norm ball admits a unital completely positive Clifford map whose compression is exactly $T_q$. For a Stinespring dilation $V^*\Gamma(x)V=T_q(x)$, Netzer's spectrahedral pencil $L$ then satisfies an explicit Gram identity
% \[
% H_{p_q,e}((0,x))[u]=\Rea\bigl(W(x)^*L(u)W(x)\bigr),
% \]
% where $H$ is the parameterized Hermite matrix and $W(x)$ is polynomial. Since $L(u)\succeq0$, realification gives a matrix SOS. Hermite--B\'ezoutian congruence, translation along $e$, and linear pullback then complete the proof. This identifies the matrix SOS certificate directly within Netzer's Clifford construction.

The third question concerns the strength of scalar SOS information. Our cubic theorem uses the fixed-direction condition, whereas the separation is first constructed in the two-direction framework of \cite{NgSaunderson2026}, where both B\'ezoutian directions range over the closed cone. For directions $u,v$ in a hyperbolicity cone, write
\[
W_{u,v}(p)=D_up\,D_vp-pD_uD_vp.
\]
This is the $(0,0)$-entry of the B\'ezoutian $B_{p,u,v}$. A matrix SOS decomposition therefore forces $W_{u,v}(p)$ to be an SOS. The issue is whether requiring these scalar SOS conditions throughout the cone recovers the matrix condition. A non-SOS Wronskian, such as the six-variable cubic example above, violates the scalar condition itself and hence does not settle this implication. To separate the two notions, the scalar conditions must all hold while a matrix obstruction survives. Our third result provides such an example.

\begin{theorem}[Two-direction weak--strong separation]\label{thm:intro-sep}
There exists a rational homogeneous polynomial
\[
P\in\Q[t,z_1,\ldots,z_{42}]
\]
of degree $46$ and a rational hyperbolicity direction $\bar e\in\Q^{43}$ such that
\[
W_{u,v}(P)\text{ is an SOS for all }u,v\in\Lam(P,\bar e),
\]
but $B_{P,u,v}$ is not a matrix SOS for some $u,v\in\Lam(P,\bar e)$.
Thus $P$ is two-direction weakly SOS-hyperbolic but not two-direction SOS-hyperbolic.
\end{theorem}

The proof separates a general localization step, stated in \cref{prop:localization}, from an explicit seed cubic built from the icosahedral graph, in the same general graph-theoretic family that underlies the high-dimensional cubic examples in \cite{Saunderson2019}. Its central Wronskian is an interior SOS quartic, certified by an explicitly reconstructible positive-definite rational Gram matrix, while an exact rational trace certificate proves that a scalar compression of its B\'ezoutian is not an SOS. A small rational polyhedral subcone is then chosen so that every Wronskian on its generators is an SOS while the B\'ezoutian obstruction persists. Multiplying by the linear forms cutting out that cone preserves weak SOS-hyperbolicity and transfers the matrix obstruction to the final degree-$46$ product. Thus the scalar two-direction theory of \cite{NgSaunderson2026} and the matrix-valued B\'ezoutian theory separate even when all data are rational.

% The linear-factor transfer is related to Saunderson's observation that multiplication by a linear form can preserve a matrix SOS obstruction \cite[Proposition~4.10]{Saunderson2019}. The localization step adds control of all scalar Wronskians on the restricted cone. Its structural content is therefore independent of the icosahedral graph: any seed satisfying the two central hypotheses can be treated by the same construction. The graph supplies explicit rational data for those hypotheses, rather than a restriction on the applicability of the localization principle.

The initial obstructing pair has the first direction $e=r_0\in\partial\Lam(P,\bar e)$, with $P(e)=0$. Although $e$ is not a hyperbolicity direction for $P$, the obstruction persists after perturbing both directions into the interior:
\[
e_\delta=e+\delta\bar e,\qquad v_\delta=r_1+\delta\bar e.
\]
For all sufficiently small $\delta>0$, \cref{cor:fixed-direction-separation} shows that $P$ is weakly SOS-hyperbolic but not SOS-hyperbolic with respect to the hyperbolicity direction $e_\delta$. Rational $\delta$ gives rational directions. This disproves unconditional equivalence of Saunderson's fixed-direction notions in Question~4.6 of \cite{Saunderson2019}.

The paper is organized as follows. \Cref{sec:prelim} fixes the B\'ezoutian, SOS, hyperzout, and SOS-hyperbolicity conventions. \Cref{sec:universal} proves the unrestricted hyperzout theorem. \Cref{sec:cubic} proves the five-variable cubic theorem. \Cref{sec:separation} gives the weak--strong separation, including the exact rational certificates. \Cref{sec:discussion} records the scope of the arguments and several remaining questions. The large Gram certificate used in the separation theorem is specified in \cref{app:certificates} in a form suitable for exact independent verification.

\section{Hyperbolicity, B\'ezoutians, and SOS certificates}\label{sec:prelim}

We fix the conventions used in all three constructions. After recalling hyperbolicity and parameterized B\'ezoutians, we distinguish the fixed-direction and two-direction SOS conditions and specify the degree profile for hyperzout representations.

\subsection{Hyperbolicity and hyperbolicity cones}

Let $p\in\R[z_1,\ldots,z_n]_d$ be homogeneous of degree $d$ and let $e\in\R^n$. We say that $p$ is \emph{hyperbolic with respect to $e$} if $p(e)>0$ and the univariate polynomial
\[
t\longmapsto p(x+te)
\]
has only real zeros for every $x\in\R^n$. The closed hyperbolicity cone $\Lam(p,e)$ is the closure of the connected component of $\{x:p(x)>0\}$ containing $e$. We write $D_up$ for the directional derivative of $p$ in the direction $u$.

Two standard facts will be used throughout. First, every interior point of $\Lam(p,e)$ is again a hyperbolicity direction and has the same hyperbolicity cone. Second, products of hyperbolic polynomials with a common hyperbolicity direction are hyperbolic, and the cone of the product is the intersection of the corresponding cones. These are classical consequences of G\aa rding's theory \cite{Garding1959}.

\subsection{B\'ezoutians}

For real univariate polynomials $a,b$ of degree at most $d$, define the $d\times d$ B\'ezoutian $B_d(a,b)$ by
\begin{equation}\label{eq:bezdef}
\frac{a(t)b(s)-b(t)a(s)}{t-s}
=\sum_{i,j=0}^{d-1}[B_d(a,b)]_{ij}t^is^j.
\end{equation}
The matrix is symmetric. Given a homogeneous polynomial $p$ and directions $u,v\in\R^n$, we set
\begin{equation}\label{eq:parambez}
B_{p,u,v}(x)\coloneqq B_d\bigl(p(x+tu),D_vp(x+tu)\bigr).
\end{equation}
When the first direction is a fixed hyperbolicity direction $e$, we also write $B_{p,e}(x)[v]$.

The upper-left entry of \eqref{eq:parambez} is
\begin{equation}\label{eq:wronskian}
[B_{p,u,v}(x)]_{00}
=D_up(x)D_vp(x)-p(x)D_uD_vp(x)
\eqqcolon W_{u,v}(p)(x).
\end{equation}
For $u,v\in\Lam(p,e)$, the parameterized B\'ezoutian is pointwise positive semidefinite (PSD); see Saunderson \cite[Section~3]{Saunderson2019} and Ng--Saunderson \cite[Section~3]{NgSaunderson2026}.

\begin{definition}\label{def:matrixsos}
A real symmetric polynomial matrix $M(x)$ is a \emph{matrix SOS} if there is a real polynomial matrix $R(x)$ such that
\[
M(x)=R(x)^\intercal R(x).
\]
A scalar polynomial is an SOS if it is a sum of squares of real polynomials.
\end{definition}

For a fixed hyperbolicity direction $e$, we use Saunderson's terminology: $p$ is \emph{SOS-hyperbolic with respect to $e$} if $B_{p,e,u}$ is a matrix SOS for every $u\in\Lam(p,e)$. It is \emph{weakly SOS-hyperbolic with respect to $e$} if $W_{e,u}(p)$ is an SOS for every such $u$. By \eqref{eq:wronskian}, SOS-hyperbolicity implies weak SOS-hyperbolicity.

For \cref{sec:separation} we also use the two-direction formulation of Ng--Saunderson: $p$ is \emph{two-direction SOS-hyperbolic} if $B_{p,u,v}$ is a matrix SOS for all $u,v\in\Lam(p,e)$, and \emph{two-direction weakly SOS-hyperbolic} if $W_{u,v}(p)$ is an SOS for all such pairs. We retain the qualifier ``two-direction'' in the separation statements. Both directions there range over the closed cone, including its boundary; neither is required to be a hyperbolicity direction. The fixed-direction notion instead requires its first direction to be the designated hyperbolicity direction. Two-direction SOS-hyperbolicity implies fixed-direction SOS-hyperbolicity, but failure of the former at a boundary first direction does not establish failure of the latter at the designated direction. Under the two-direction weak condition, the perturbation in \cref{cor:fixed-direction-separation} yields a fixed-direction separation at a nearby interior first direction.

\subsection{Hyperzout representations}\label{sec:prelim-hyperzout}

Let $\mathcal F_{m,j}$ be the vector space of real homogeneous forms of degree $j$ in $m$ variables, with the zero polynomial included in every $\mathcal F_{m,j}$. For positive integers $d,k,\mu$ with $\mu\le d$, set
\begin{equation}\label{eq:Tspace}
\calT^{m,d}_{\mu,k}
\coloneqq\{0\}^{d-\mu}\times \mathcal F_{m,0}\times\mathcal F_{m,k}
\times\cdots\times\mathcal F_{m,(\mu-1)k}.
\end{equation}
Let
\[
\calP_{m,2r}\coloneqq\{f\in\mathcal F_{m,2r}:f(x)\ge0\text{ for all }x\in\R^m\}.
\]

\begin{definition}[Ng--Saunderson]\label{def:hyperzout}
A form $f\in\mathcal F_{m,2r}$ is \emph{hyperzout} if there exist positive integers $n,d,k,\mu$, a hyperbolic form $p\in\R[z_1,\ldots,z_n]_d$ with hyperbolicity direction $e$, directions $u,v\in\Lam(p,e)$, a homogeneous polynomial map $\phi\in(\mathcal F_{m,k})^n$, and a vector $\xi\in\calT^{m,d}_{\mu,k}$ such that
\[
r=(\mu-1)k,
\qquad
f(x)=\xi(x)^\intercal B_{p,u,v}(\phi(x))\xi(x).
\]
We denote the set of such forms by $\calZ_{m,2r}$.
\end{definition}

The degree-restricted subclass considered in \cite{NgSaunderson2026} additionally requires
\begin{equation}\label{eq:restricted-degree}
\mu=2\qquad\text{or}\qquad d\le 2\mu-3.
\end{equation}
We denote this subclass by $\calB_{m,2r}$. Since $B_{p,u,v}$ is pointwise PSD on real inputs whenever $u,v$ lie in the hyperbolicity cone, every hyperzout form is nonnegative:
\begin{equation}\label{eq:ZsubsetP}
\calZ_{m,2r}\subseteq\calP_{m,2r}.
\end{equation}

\section{Universal hyperzout representations}\label{sec:universal}

To prove \cref{thm:intro-univ}, we first construct an auxiliary hyperbolic polynomial by perturbing a square. A root-moment factorization then recovers the target form as a B\'ezoutian compression, and a triangular coefficient calculation gives the required degree profile.

\subsection{A homogeneous real-rooted perturbation}

The auxiliary polynomial is obtained from a real-rooted univariate polynomial whose roots remain real under small level perturbations.

\begin{lemma}\label{lem:level-perturb}
Let $H\in\R[z]$ be monic of degree $r$ and suppose that $H$ has $r$ distinct real roots. Then there exists $\tau>0$ such that $H(z)-c$ has $r$ distinct real roots for every $|c|<\tau$.
\end{lemma}

\begin{proof}
Write the roots as $\alpha_1<\cdots<\alpha_r$. Since they are simple, for each $i$ one may choose an interval $(a_i,b_i)$ containing $\alpha_i$, with the intervals pairwise disjoint, such that $H(a_i)H(b_i)<0$. Choose
\[
0<\tau<\min_i\{|H(a_i)|,|H(b_i)|\}.
\]
If $|c|<\tau$, then $H(a_i)-c$ and $H(b_i)-c$ still have opposite signs. The intermediate value theorem gives one zero of $H-c$ in each $(a_i,b_i)$. There are $r$ such intervals and the degree is $r$, so all roots are real and simple.
\end{proof}

To apply this stability observation homogeneously, fix $r\ge1$ and write
\[
s=\lfloor r/2\rfloor,\qquad \delta=r-2s\in\{0,1\},
\]
and define
\begin{equation}\label{eq:Hr}
H_r(z)\coloneqq z^\delta\prod_{j=1}^s(z^2-j^2).
\end{equation}
Its roots are $\pm1,\ldots,\pm s$, together with $0$ when $r$ is odd. Let $\tau_r>0$ be supplied by \cref{lem:level-perturb}. For $x=(x_1,\ldots,x_m)$ set
\begin{equation}\label{eq:hr}
b(x)\coloneqq\sum_{i=1}^m x_i^2,
\qquad
h(t,x)\coloneqq t^\delta\prod_{j=1}^s\bigl(t^2-j^2b(x)\bigr).
\end{equation}
The resulting $h$ is homogeneous of total degree $r$ and monic in $t$.

\begin{proposition}\label{prop:aux-hyperbolic}
Let $0\ne f\in\calP_{m,2r}$ and set
\[
M\coloneqq \max_{\|x\|=1} f(x)>0.
\]
Choose
\[
0<\eps<\frac{\tau_r^2}{M}
\]
and define
\begin{equation}\label{eq:paux}
p(t,x)\coloneqq h(t,x)^2-\eps f(x).
\end{equation}
Then $p$ is a homogeneous form of degree $2r$ in $m+1$ variables, hyperbolic with respect to $e=(1,0,\ldots,0)$, and $p(e)=1$.
\end{proposition}

\begin{proof}
Compactness of the unit sphere gives a finite maximum $M$, and $M>0$ because $f$ is nonzero and nonnegative. Let $x\ne0$, put $\rho=\|x\|$, and use homogeneity to write
\[
h(t,x)=\rho^r H_r(t/\rho),
\qquad
f(x)=\rho^{2r}f(x/\rho).
\]
If
\[
\alpha\coloneqq \sqrt{\eps f(x/\rho)},
\]
then $0\le\alpha\le\sqrt{\eps M}<\tau_r$, and
\[
p(t,x)=\rho^{2r}\bigl(H_r(t/\rho)-\alpha\bigr)
\bigl(H_r(t/\rho)+\alpha\bigr).
\]
Both factors have only real roots by \cref{lem:level-perturb}; when $\alpha=0$ the two factors coincide, but the roots remain real with multiplicity. If $x=0$, then $p(t,0)=t^{2r}$. Thus every fiber $t\mapsto p(t,x)$ is real-rooted.

For an arbitrary base point $X=(s,x)$,
\[
p(\lambda e-X)=p(\lambda-s,-x).
\]
The preceding real-rootedness is unaffected by the translation $\lambda\mapsto\lambda-s$. Since $p(e)=1$, this is precisely hyperbolicity with respect to $e$.
\end{proof}

\subsection{A B\'ezoutian moment factorization}

The next step relates the B\'ezoutian to moments of its roots. We use the coefficient vector
\[
V_D(t)\coloneqq (1,t,\ldots,t^{D-1})^\intercal .
\]

\begin{lemma}\label{lem:root-moment}
Let $P(t)=\sum_{j=0}^D a_jt^j$ be a monic real-rooted polynomial of degree $D$, and list its roots with multiplicity as $\lambda_1,\ldots,\lambda_D$. Let
\[
C\coloneqq B_D(P,1),
\qquad
H\coloneqq \sum_{\ell=1}^D V_D(\lambda_\ell)V_D(\lambda_\ell)^\intercal .
\]
Then
\begin{equation}\label{eq:moment-factor}
B_D(P,P')=CHC.
\end{equation}
In particular, $B_D(P,P')\succeq0$.
\end{lemma}

\begin{proof}
For each listed root, including repetitions, set $Q_\ell(t)=P(t)/(t-\lambda_\ell)$. The product rule gives
\[
P'(t)=\sum_{\ell=1}^D Q_\ell(t).
\]
Consequently,
\begin{align*}
\frac{P(t)P'(s)-P(s)P'(t)}{t-s}
&=\sum_{\ell=1}^D
\frac{P(t)Q_\ell(s)-P(s)Q_\ell(t)}{t-s}\\
&=\sum_{\ell=1}^D Q_\ell(t)Q_\ell(s).
\end{align*}
We express the quotient polynomials in coefficient form using the definition of $C$:
\[
\frac{P(t)-P(s)}{t-s}=V_D(t)^\intercal CV_D(s).
\]
In this identity, set $s=\lambda_\ell$. Since $P(\lambda_\ell)=0$, its left-hand side becomes $P(t)/(t-\lambda_\ell)=Q_\ell(t)$, and hence
\[
Q_\ell(t)=V_D(t)^\intercal CV_D(\lambda_\ell).
\]
Substituting this expression for each $Q_\ell$ into the preceding sum gives
\begin{align*}
\frac{P(t)P'(s)-P(s)P'(t)}{t-s}
&=\sum_{\ell=1}^D
\bigl(V_D(t)^\intercal CV_D(\lambda_\ell)\bigr)
\bigl(V_D(s)^\intercal CV_D(\lambda_\ell)\bigr)\\
&=V_D(t)^\intercal  C
\left(\sum_{\ell=1}^D V_D(\lambda_\ell)V_D(\lambda_\ell)^\intercal \right)
C V_D(s)\\
&=V_D(t)^\intercal CHCV_D(s),
\end{align*}
Here we use the definition of $H$ and the symmetry $C=C^\intercal $. Comparing coefficients of $t^is^j$ gives $B_D(P,P')=CHC$, proving \eqref{eq:moment-factor}. Multiplicities are included in both the product-rule sum and $H$. Only the polynomial quotients by $t-\lambda_\ell$ are used to define $Q_\ell$; no division by $P'(\lambda_\ell)$ is required, so the argument also covers repeated roots.
\end{proof}

To obtain a hyperzout representation, the compression vector must also satisfy the prescribed degree profile. The following triangular calculation supplies this condition.

\begin{lemma}\label{lem:degree-profile}
Let $D=2r$ and
\[
P(t,x)=\sum_{j=0}^D a_j(x)t^j
\]
be homogeneous of total degree $D$ and monic in $t$. Let
\[
h(t,x)=\sum_{j=0}^r z_j(x)t^j
\]
be homogeneous of total degree $r$ and monic in $t$. Pad its coefficient vector to length $D$:
\[
z(x)=(z_0,\ldots,z_r,0,\ldots,0)^\intercal ,
\qquad
C(x)\coloneqq B_D(P(\cdot,x),1).
\]
Then $C(x)$ is invertible and $C(x)^{-1}$ has polynomial entries, and
\begin{equation}\label{eq:eta-profile}
\eta(x)\coloneqq C(x)^{-1}z(x)
\in \{0\}^{r-1}\times\mathcal F_{m,0}\times\mathcal F_{m,1}
\times\cdots\times\mathcal F_{m,r}.
\end{equation}
Moreover, the first nonzero-profile entry is $\eta_{r-1}=1$ when indices start at $0$.
\end{lemma}

\begin{proof}
From \eqref{eq:bezdef}, with rows and columns indexed from $0$,
\begin{equation}\label{eq:Centries}
C_{ij}=\begin{cases}
a_{i+j+1},&i+j\le D-1,\\
0,&i+j>D-1.
\end{cases}
\end{equation}
Let $J$ reverse the order of the columns. Then $CJ$ is upper triangular with all diagonal entries equal to $a_D=1$. Hence
\[
\det C=(-1)^{D(D-1)/2},
\]
a nonzero constant. The adjugate formula therefore shows that $C^{-1}$ is polynomial.

It remains to determine the degrees. Homogeneity gives
\[
a_j\in\mathcal F_{m,D-j},
\qquad
z_j\in\mathcal F_{m,r-j}.
\]
Extend $z_j=0$ for $j>r$. The $i$th row of $C\eta=z$ is
\begin{equation}\label{eq:Ceta-row}
\sum_{j=0}^{D-1-i}a_{i+j+1}\eta_j=z_i.
\end{equation}
Starting at $i=D-1$ and descending to $i=r+1$ yields
\[
\eta_0=\cdots=\eta_{r-2}=0.
\]
At $i=r$, \eqref{eq:Ceta-row} gives $\eta_{r-1}=z_r=1$. For $q=1,\ldots,r$, set $i=r-q$. Then
\[
\eta_{r-1+q}
=z_{r-q}-\sum_{j=r-1}^{r+q-2}a_{r-q+j+1}\eta_j.
\]
Assuming inductively that $\eta_j\in\mathcal F_{m,j-r+1}$, each nonzero product in the sum has degree
\[
\bigl(2r-(r-q+j+1)\bigr)+(j-r+1)=q,
\]
and $z_{r-q}$ is either zero or homogeneous of degree $q$. Hence $\eta_{r-1+q}$ is zero or homogeneous of degree $q$. Since zero belongs to $\mathcal F_{m,q}$, this proves $\eta_{r-1+q}\in\mathcal F_{m,q}$ and establishes \eqref{eq:eta-profile}.
\end{proof}

\begin{remark}\label{rem:degree-profile-parity}
 Equation~\eqref{eq:hr} gives $h(-t,x)=(-1)^r h(t,x)$, so $p(t,x)=h(t,x)^2-\eps f(x)$ is even in $t$. Thus, in the notation of \cref{lem:degree-profile}, $a_j(x)=0$ for every odd $j$. Accordingly, $C_{ij}=0$ when $i$ and $j$ have the same parity. Since $z_i=0$ unless $i\equiv r\pmod 2$, solving $C\eta=z$ by parity blocks also gives
\[
\eta_j(x)=0\qquad\text{whenever }j\equiv r\pmod 2.
\]
The general degree-profile statement and proof above do not require this specialization.
\end{remark}

\subsection{Universality}

We now combine the auxiliary polynomial, the moment factorization, and the degree calculation to prove the first main theorem.

\begin{theorem}\label{thm:hyperzout-universal}
For all $m,r\ge1$,
\[
\calZ_{m,2r}=\calP_{m,2r}.
\]
For every nonzero $f\in\calP_{m,2r}$ one may take
\begin{equation}\label{eq:univ-params}
n=m+1,\quad d=2r,\quad k=1,\quad\mu=r+1,
\quad e=u=v=(1,0,\ldots,0),\quad \phi(x)=(0,x).
\end{equation}
\end{theorem}

\begin{proof}
The inclusion $\calZ_{m,2r}\subseteq\calP_{m,2r}$ is \eqref{eq:ZsubsetP}. For the converse, first suppose $f\ne0$. Choose $\eps$ and $p$ from \cref{prop:aux-hyperbolic}, and set $D=2r$. For fixed $x$, write
\[
p_x(t)\coloneqq p(t,x),
\qquad
C(x)\coloneqq B_D(p_x,1).
\]
Let $\lambda_1(x),\ldots,\lambda_D(x)$ be the roots of $p_x$, counted with multiplicity, and set
\[
H(x)\coloneqq \sum_{\ell=1}^D V_D(\lambda_\ell(x))V_D(\lambda_\ell(x))^\intercal .
\]
Only the pointwise root values are used, so no regularity in $x$ is required. The factorization in \cref{lem:root-moment} gives
\begin{equation}\label{eq:BPPC}
B_D(p_x,p_x')=C(x)H(x)C(x).
\end{equation}

Let $z(x)$ be the coefficient vector of $h(t,x)$, padded to length $D$. Then
\[
z(x)^\intercal V_D(t)=h(t,x).
\]
Since every root of $p_x$ satisfies $h(\lambda_\ell(x),x)^2=\eps f(x)$,
\begin{equation}\label{eq:zHz}
z(x)^\intercal H(x)z(x)
=\sum_{\ell=1}^D h(\lambda_\ell(x),x)^2
=D\eps f(x).
\end{equation}
Define
\begin{equation}\label{eq:xi-univ}
\xi(x)\coloneqq \frac{1}{\sqrt{D\eps}}C(x)^{-1}z(x).
\end{equation}
By \cref{lem:degree-profile}, $\xi$ is polynomial and
\[
\xi\in\calT^{m,2r}_{r+1,1}.
\]
Combining \eqref{eq:BPPC} and \eqref{eq:zHz} gives
\[
\xi(x)^\intercal B_D(p_x,p_x')\xi(x)=f(x).
\]
Finally,
\[
p(\phi(x)+te)=p_x(t),
\qquad
D_ep(\phi(x)+te)=p_x'(t),
\]
so the displayed identity is exactly the hyperzout representation required by \cref{def:hyperzout} with the parameters \eqref{eq:univ-params}.

If $f=0$, take $p=h^2$ and $\xi=0$. The form $p$ is hyperbolic because every fiber of $h$ is real-rooted, and the zero vector belongs to the required degree-profile space. This completes the proof.
\end{proof}

\begin{corollary}\label{cor:Zclosed}
For fixed $m$ and $r$, the unrestricted hyperzout set $\calZ_{m,2r}$ is a closed convex cone.
\end{corollary}

\begin{proof}
By \cref{thm:hyperzout-universal}, it equals the cone $\calP_{m,2r}$ of nonnegative forms, which is the intersection of the closed half-spaces $\{f:f(x)\ge0\}$ indexed by $x\in\R^m$.
\end{proof}

\begin{remark}\label{rem:degree-boundary}
Here ``unrestricted'' describes the hyperzout class, not a degree-free certificate. The construction in \cref{thm:hyperzout-universal} has
\[
d=2r,\qquad k=1,\qquad \mu=r+1,
\qquad\text{so}\qquad dk=2r.
\]
When $r\ge2$, one has $\mu=r+1\ge3$, so neither branch of \eqref{eq:restricted-degree} holds. Thus the construction lies exactly one auxiliary degree above the inequality $d\le2\mu-3$. In the notation of \cite{NgSaunderson2026}, their dimension-counting obstruction applies to $\calB_{m,2r}$, not to the unrestricted $\calZ_{m,2r}$. The theorem does not assert that $2r$ is the smallest possible auxiliary degree for every nonnegative form.
\end{remark}

\begin{remark}\label{rem:nonalgorithmic}
The construction presupposes $f\ge0$. This assumption is used to make $\sqrt{\eps f(x/\|x\|)}$ real and to ensure hyperbolicity of $p=h^2-\eps f$. If $f(x_0)<0$ at some point, then $h(t,x_0)^2-\eps f(x_0)>0$ for every real $t$, so that fiber has no real roots and the same auxiliary construction cannot be hyperbolic. Thus \cref{thm:hyperzout-universal} provides an existential representation rather than an efficient nonnegativity test.
\end{remark}

\begin{example}\label{ex:quartic-hyperzout}
Let $0\ne f\in\calP_{m,4}$, put $b=\|x\|^2$, and choose $0<\eps<1/M$ with $M=\max_{\|x\|=1}f(x)$. Then
\[
p(t,x)=(t^2-b)^2-\eps f(x)
\]
is hyperbolic with respect to $e=(1,0,\ldots,0)$. Writing $a=\eps f(x)$, one finds
\[
B_4(p_x,p_x')=4
\begin{pmatrix}
 b^3-ab&0&a-b^2&0\\
 0&a+b^2&0&-b\\
 a-b^2&0&b&0\\
 0&-b&0&1
\end{pmatrix}.
\]
The vector
\[
\xi(x)=\frac1{2\sqrt\eps}(0,1,0,b)^\intercal 
\in \{0\}\times\mathcal F_{m,0}\times\mathcal F_{m,1}\times\mathcal F_{m,2}
\]
satisfies $\xi^\intercal B_4(p_x,p_x')\xi=f(x)$.
\end{example}

\section{Five-variable hyperbolic cubics are SOS-hyperbolic}\label{sec:cubic}

We prove \cref{thm:intro-cubic} by converting Netzer's canonical matrix compression \cite{Netzer2026} into a polynomial Hermite Gram factorization. We first record the compression property, then establish the Gram identity, and finally transfer the resulting SOS decomposition to arbitrary coordinates and hyperbolicity directions.

\subsection{Normal form and the canonical pencil}

After normalizing $p(e)=1$, applying a linear change of variables, and quotienting transverse directions on which $p$ is constant, Netzer's normal-form result gives either the pure cube $t^3$ or
\begin{equation}\label{eq:cubic-normal}
p_q(t,x)=t^3-3t\|x\|^2+2q(x),
\qquad x\in\R^m,
\end{equation}
with $1\le m\le4$ and
\begin{equation}\label{eq:Hm}
q\in\calH_m\coloneqq \{q\in\R[x_1,\ldots,x_m]_3:\|q\|_\infty\le1\},
\qquad
\|q\|_\infty\coloneqq \max_{\|x\|=1}|q(x)|.
\end{equation}
Conversely, \eqref{eq:cubic-normal} is hyperbolic with respect to $e_0=(1,0)$ if and only if $q\in\calH_m$ \cite[Proposition~1.1]{Netzer2026}.

Let $B_q$ be the symmetric trilinear form associated with $q$, so $q(x)=B_q(x,x,x)$. Define the canonical symmetric pencil $T_q$ by
\begin{equation}\label{eq:Tq-def}
u^\intercal T_q(x)v=B_q(x,u,v).
\end{equation}
Equivalently,
\begin{equation}\label{eq:Tq-hess}
T_q(x)=\frac16\nabla^2q(x),
\qquad
q(x)=x^\intercal T_q(x)x,
\qquad
\nabla q(x)=3T_q(x)x.
\end{equation}
Write $T_q(x)=\sum_i x_iT_{q,i}$.

\subsection{Canonical Clifford compressions}

Let $\calC_m$ be the universal complex Clifford $C^*$-algebra generated by self-adjoint elements $c_1,\ldots,c_m$ satisfying
\begin{equation}\label{eq:Clifford}
c_ic_j+c_jc_i=2\delta_{ij}I.
\end{equation}
Set $c(x)=\sum_i x_ic_i$. A Clifford realization in Netzer's sense is a unital completely positive map $\Theta:\calC_m\to M_m(\C)$ such that
\[
q(x)=x^*\Theta(c(x))x.
\]
A Stinespring dilation \cite{Stinespring1955} gives a finite-dimensional Hilbert space $K$, a $*$-representation $\pi:\calC_m\to B(K)$, and an isometry $V:\C^m\to K$ such that
\[
\Theta(a)=V^*\pi(a)V.
\]
Writing $\Gamma_i=\pi(c_i)$ and $\Gamma(x)=\sum_i x_i\Gamma_i$, we have
\begin{equation}\label{eq:Gamma-rel}
\Gamma(x)^2=\|x\|^2I,
\qquad
\Gamma(x)\Gamma(y)+\Gamma(y)\Gamma(x)=2x^\intercal yI.
\end{equation}

The Gram factorization requires more than the scalar identity $x^*V^*\Gamma(x)Vx=q(x)$. Specifically, it uses
\begin{equation}\label{eq:canonical-compression}
V^*\Gamma(x)V=T_q(x).
\end{equation}
Netzer's Proposition~3.8 gives this identity, and the argument of his Theorem~3.9 extends it to the full norm ball. We recall the construction to identify precisely the input needed below.

For $q\in\calH_4$, let
\[
Z(q)\coloneqq \{v\in\Sph^3:q(v)=1\}.
\]
A five-tuple in $Z(q)$ is \emph{poised} if every four of its members are linearly independent. Following Netzer, $q$ is a \emph{five-contact cubic} if $Z(q)$ contains a poised five-tuple. Five-contact cubics are dense in $\Ext(\calH_4)$ \cite[Proposition~3.3]{Netzer2026}.

\begin{lemma}[Netzer, Proposition~3.8]\label{prop:five-contact}
Let $q\in\calH_4$ be five-contact. Then there exists a finite-dimensional complex Hilbert space $K$, an isometry $V:\C^4\to K$, and a Clifford representation $x\mapsto\Gamma(x)$ such that \eqref{eq:canonical-compression} holds.
\end{lemma}

\begin{proof}
After an orthogonal change of coordinates, one contact point may be placed at $e_1$. Writing $\R^4=\R e_1\oplus U$ and $x=(s,y)$, the contact condition yields
\begin{equation}\label{eq:q-adapted}
q(s,y)=s^3+3s\,y^\intercal Dy+r(y),
\end{equation}
where $D=D^\intercal $, $r$ is cubic on $U$, and $T(y)=\frac16\nabla^2r(y)$. Let
\[
P=\frac{I+D}{2},\qquad Q=\frac{I-D}{2}.
\]
The quaternionic determinantal construction in \cite[Propositions~3.4 and 3.6]{Netzer2026} gives a finite-dimensional real Hilbert space $F$, a unit vector $\xi\in F$, and maps $a:U\to F$, $J:U\to\operatorname{Hom}(U,F)$ satisfying
\begin{gather}\label{eq:netzer-identities1}
a^*\xi=0,\qquad J(y)^*\xi=Dy,\qquad a^*a=P,\\
\label{eq:netzer-identities2}
a^*J(y)+J(y)^*a=T(y),\qquad J(y)^*J(y)=\|y\|^2Q.
\end{gather}
In the five-contact case, $Q\succ0$. Netzer's Clifford-extension lemma provides Hilbert spaces $K_-,K_+$, isometric embeddings of $U_\C$ and $F_\C$, and a linear map $C(y):K_-\to K_+$ with
\[
C(y)^*C(y)=\|y\|^2I,
\qquad
C(y)C(y)^*=\|y\|^2I,
\]
and, after identifying $F_\C$ with its image in $K_+$,
\[
C(y)b(u)=J(y)u,
\qquad
b\coloneqq \iota_-Q^{1/2}.
\]
On $K=K_+\oplus K_-$ set
\[
\Gamma_0=\begin{pmatrix}I&0\\0&-I\end{pmatrix},
\qquad
\Gamma(y)=\begin{pmatrix}0&C(y)\\C(y)^*&0\end{pmatrix}.
\]
These matrices satisfy the Clifford relations. Define
\[
V(\alpha,u)=\binom{\alpha\xi+a(u)}{b(u)}.
\]
Because $a^*a+b^*b=P+Q=I$ and $a^*\xi=0$, the map $V$ is an isometry. Using \eqref{eq:netzer-identities1}--\eqref{eq:netzer-identities2} and $C(y)b(u)=J(y)u$, block multiplication gives
\[
V^*\bigl(s\Gamma_0+\Gamma(y)\bigr)V
=
\begin{pmatrix}
s&(Dy)^\intercal \\
Dy&sD+T(y)
\end{pmatrix}.
\]
By the Hessian formula for \eqref{eq:q-adapted}, the right-hand side is $T_q(s,y)$.
\end{proof}

To extend this construction beyond five-contact cubics, define
\begin{equation}\label{eq:E4}
\calE_4\coloneqq \left\{q\in\R[x_1,\ldots,x_4]_3:
\begin{array}{l}
\text{there is a unital completely positive (UCP) map }\Theta:\calC_4\to M_4(\C)\\
\text{with }\Theta(c_i)=T_{q,i}\text{ for }i=1,\ldots,4
\end{array}
\right\}.
\end{equation}

\begin{lemma}\label{lem:O4-inv}
The set $\calE_4$ is invariant under orthogonal changes of variables.
\end{lemma}

\begin{proof}
Let $U\in O(4)$ and define $q_U(x)=q(U^\intercal x)$. The Hessian transformation rule gives
\begin{equation}\label{eq:Ttransform}
T_{q_U}(x)=UT_q(U^\intercal x)U^\intercal .
\end{equation}
Let $\alpha_{U^\intercal }$ be the Clifford automorphism determined by
\[
\alpha_{U^\intercal }(c(x))=c(U^\intercal x).
\]
If $\Theta$ witnesses $q\in\calE_4$, then
\[
\Theta_U(a)\coloneqq U\Theta(\alpha_{U^\intercal }(a))U^\intercal 
\]
is again UCP, and \eqref{eq:Ttransform} gives $\Theta_U(c(x))=T_{q_U}(x)$.
\end{proof}

\begin{lemma}\label{thm:E4H4}
One has \[
\calE_4=\calH_4.
\]
Equivalently, every $q\in\calH_4$ admits a Stinespring realization satisfying \eqref{eq:canonical-compression} and \eqref{eq:Gamma-rel}.
\end{lemma}

\begin{proof}
We follow \cite[Proposition~3.8 and Theorem~3.9]{Netzer2026}, retaining the canonical compression under convex combinations and limits. Convexity of $\calE_4$ follows because convex combinations of UCP maps are UCP and $q\mapsto T_q$ is linear. To prove closedness, let $q_j\to q$ with $q_j\in\calE_4$, and choose UCP maps $\Theta_j$ realizing these compressions. The finite-dimensional set of UCP maps $\calC_4\to M_4(\C)$ is compact, as follows from the Choi-matrix description on the matrix summands of $\calC_4$. Passing to a subsequence gives $\Theta_j\to\Theta$, and continuity yields $\Theta(c_i)=T_{q,i}$. Thus $\calE_4$ is closed.

By \cref{prop:five-contact,lem:O4-inv}, every five-contact cubic belongs to $\calE_4$. Since these cubics are dense in $\Ext(\calH_4)$, closedness gives
\[
\Ext(\calH_4)\subseteq\calE_4.
\]
The set $\calH_4$ is compact and convex, hence it is the convex hull of its extreme points. Convexity of $\calE_4$ therefore gives $\calH_4\subseteq\calE_4$.

Conversely, if $q\in\calE_4$, then a UCP map is contractive on self-adjoint elements, so
\[
\|T_q(x)\|=\|\Theta(c(x))\|\le\|c(x)\|=\|x\|.
\]
Using $q(x)=x^\intercal T_q(x)x$ gives $|q(x)|\le\|x\|^3$, hence $q\in\calH_4$. Thus $\calE_4=\calH_4$.
\end{proof}

If $\mathcal R_4$ denotes Netzer's class defined by the scalar Clifford identity, then $\calE_4\subseteq\mathcal R_4$. Although $\calE_4=\calH_4$ is formally stronger than the stated equality $\mathcal R_4=\calH_4$, it follows from Netzer's proof. Our next step is to derive the Hermite Gram identity from this compression property.

\subsection{Netzer's pencil and the Hermite Gram identity}

With the compression property established, fix $q\in\calH_4$ and a canonical Stinespring realization. Netzer's spectrahedral pencil is
\begin{equation}\label{eq:Netzer-pencil}
L_q(t,x)=
\begin{pmatrix}
\frac{t}{3}&\sqrt{\frac23}(Vx)^*\\[1mm]
\sqrt{\frac23}Vx&tI+\Gamma(x)
\end{pmatrix}.
\end{equation}
Its positive semidefinite locus is the hyperbolicity cone:
\begin{equation}\label{eq:cone-pencil}
\Lam(p_q,e_0)=\{(t,x):L_q(t,x)\succeq0\}.
\end{equation}
For a fixed cone point $u=(a,b)$, write
\begin{equation}\label{eq:L-u}
L(u)=
\begin{pmatrix}
\frac{a}{3}&\sqrt{\frac23}B^*\\[1mm]
\sqrt{\frac23}B&aI+J
\end{pmatrix}\succeq0,
\qquad B\coloneqq Vb,\quad J\coloneqq \Gamma(b).
\end{equation}

To compare this pencil with the Hermite matrix, for $x\in\R^4$ write
\begin{equation}\label{eq:notation-cubic}
r^2\coloneqq \|x\|^2,\qquad \beta\coloneqq b^\intercal x,\qquad \gamma\coloneqq D_bq(x),
\qquad A\coloneqq Vx,\quad G\coloneqq \Gamma(x).
\end{equation}

\begin{lemma}\label{lem:clifford-identities}
With the notation above,
\begin{gather}
A^*A=r^2,\qquad A^*B=\beta,\qquad A^*GA=q(x),\label{eq:cliff1}\\
B^*GA=\frac{\gamma}{3},\qquad A^*JA=\frac{\gamma}{3},\qquad JG+GJ=2\beta I,\label{eq:cliff2}\\
\Rea(A^*JGA)=\beta r^2,
\qquad
A^*GJGA=2\beta q(x)-\frac{r^2\gamma}{3}.\label{eq:cliff3}
\end{gather}
\end{lemma}

\begin{proof}
The first two identities in \eqref{eq:cliff1} follow from $V^*V=I$. The third follows from $V^*\Gamma(x)V=T_q(x)$ and \eqref{eq:Tq-hess}. Since
\[
\gamma=D_bq(x)=3B_q(b,x,x),
\]
symmetry of $B_q$ gives
\[
B^*GA=b^\intercal T_q(x)x=B_q(b,x,x)=\frac\gamma3,
\]
and similarly $A^*JA=x^\intercal T_q(b)x=\gamma/3$. The Clifford relation gives $JG+GJ=2\beta I$, hence
\[
2\Rea(A^*JGA)=A^*(JG+GJ)A=2\beta r^2.
\]
Finally, multiplying $JG+GJ=2\beta I$ on the left by $G$ gives
\[
GJG=2\beta G-r^2J,
\]
which yields the last identity.
\end{proof}

We now compute the Hermite matrix and compare its entries with the Clifford expressions above. For a univariate pair $a,b$ with $a$ monic of degree $d$, expand at infinity
\[
\frac{b(s)}{a(s)}=\sum_{j\ge0}\eta_j s^{-j-1}.
\]
The associated $d\times d$ Hermite matrix is $[\eta_{i+j}]_{i,j=0}^{d-1}$. For hyperbolic $p$ we denote the corresponding parameterized matrix by $H_{p,e}(x)[u]$.

On the transverse section $t=0$, set
\[
f(s)\coloneqq p_q(s,x)=s^3-3r^2s+2q(x).
\]
The directional derivative along $u=(a,b)$ is
\[
D_up_q(s,x)=3as^2-6\beta s-3ar^2+2\gamma.
\]
Writing
\[
\frac{D_up_q(s,x)}{p_q(s,x)}
=\frac{\eta_0}{s}+\frac{\eta_1}{s^2}+\frac{\eta_2}{s^3}
+\frac{\eta_3}{s^4}+\frac{\eta_4}{s^5}+O(s^{-6})
\]
and comparing coefficients gives
\begin{equation}\label{eq:etas-cubic}
\eta_0=3a,
\quad \eta_1=-6\beta,
\quad \eta_2=6ar^2+2\gamma,
\quad \eta_3=-6aq-18\beta r^2,
\quad \eta_4=18ar^4+6\gamma r^2+12\beta q.
\end{equation}
Hence
\begin{equation}\label{eq:H-cubic}
H(x;u)\coloneqq H_{p_q,e_0}((0,x))[u]
=\begin{pmatrix}
3a&-6\beta&6ar^2+2\gamma\\
-6\beta&6ar^2+2\gamma&-6aq-18\beta r^2\\
6ar^2+2\gamma&-6aq-18\beta r^2&18ar^4+6\gamma r^2+12\beta q
\end{pmatrix}.
\end{equation}

To express the Hermite matrix as a Gram matrix, define the polynomial block matrix with columns in $\C\oplus K$ by
\begin{equation}\label{eq:Wmatrix}
W(x)=
\begin{pmatrix}
3&0&6r^2\\
0&-\sqrt6\,A&\sqrt6\,GA
\end{pmatrix}.
\end{equation}

\begin{proposition}\label{prop:Gram-cubic}
For every $q\in\calH_4$, $u\in\Lam(p_q,e_0)$, and $x\in\R^4$,
\begin{equation}\label{eq:Gram-cubic}
H(x;u)=\Rea\bigl(W(x)^*L(u)W(x)\bigr).
\end{equation}
\end{proposition}

\begin{proof}
Let $w_1,w_2,w_3$ be the columns of $W$. We verify the independent entries using \cref{lem:clifford-identities}:
\begin{align*}
w_1^*L(u)w_1&=3a,\\
\Rea(w_1^*L(u)w_2)&=-6\beta,\\
\Rea(w_1^*L(u)w_3)&=6ar^2+6B^*GA=6ar^2+2\gamma,\\
w_2^*L(u)w_2&=6A^*(aI+J)A=6ar^2+2\gamma.
\end{align*}
For the $(2,3)$ entry,
\begin{align*}
\Rea(w_2^*L(u)w_3)
&=-12r^2\Rea(A^*B)-6aA^*GA-6\Rea(A^*JGA)\\
&=-6aq-18\beta r^2.
\end{align*}
Finally,
\begin{align*}
w_3^*L(u)w_3
&=12ar^4+24r^2\Rea(B^*GA)+6A^*G(aI+J)GA\\
&=12ar^4+8\gamma r^2+6ar^4
+6\left(2\beta q-\frac{r^2\gamma}{3}\right)\\
&=18ar^4+6\gamma r^2+12\beta q.
\end{align*}
These are the six independent entries of \eqref{eq:H-cubic}.
\end{proof}

\begin{corollary}\label{cor:H-matrix-sos}
For every $q\in\calH_4$ and every $u\in\Lam(p_q,e_0)$, the matrix polynomial $H_{p_q,e_0}((0,x))[u]$ is a matrix SOS.
\end{corollary}

\begin{proof}
Since $L(u)\succeq0$ is constant in $x$, choose a complex matrix $C$ with $L(u)=C^*C$. Set
\[
Z(x)=CW(x)=P(x)+iQ(x)
\]
with $P,Q$ real polynomial matrices. Then \cref{prop:Gram-cubic} gives
\[
H(x;u)=\Rea(Z^*Z)=P^\intercal P+Q^\intercal Q.
\]
\end{proof}

\subsection{From the Hermite matrix to an arbitrary cubic}

It remains to transfer the transverse Hermite factorization to the full B\'ezoutian and then to an arbitrary cubic. We begin with the coefficient matrix in the Hermite--B\'ezoutian congruence \cite[Proposition~3.4]{Saunderson2019}. Write
\[
a(\tau)=p(x+\tau e)=\sum_{j=0}^d a_j(x)\tau^j,
\qquad a_d=p(e)=1,
\]
and put $b(\tau)=D_up(x+\tau e)$. With indices starting at $0$, the symmetric coefficient matrix is
\begin{equation}\label{eq:explicit-HB-matrix}
C(x)=B_d(a,1)=
\begin{pmatrix}
a_1&a_2&\cdots&a_{d-1}&1\\
a_2&a_3&\cdots&1&0\\
\vdots&\vdots&\cdots&\vdots&\vdots\\
a_{d-1}&1&\cdots&0&0\\
1&0&\cdots&0&0
\end{pmatrix},\qquad
C_{ij}=\begin{cases}a_{i+j+1},&i+j\le d-1,\\0,&i+j>d-1.\end{cases}
\end{equation}
Its determinant is $(-1)^{d(d-1)/2}$, so it is unimodular. If $a$ has distinct roots $\rho_\ell$, partial fractions give
\[
\frac{b(\tau)}{a(\tau)}=\sum_\ell\frac{\omega_\ell}{\tau-\rho_\ell},
\qquad \omega_\ell=\frac{b(\rho_\ell)}{a'(\rho_\ell)}.
\]
Thus $H=\sum_\ell\omega_\ell V_d(\rho_\ell)V_d(\rho_\ell)^\intercal$. The coefficient vector of $a(\tau)/(\tau-\rho_\ell)$ is $CV_d(\rho_\ell)$, as in \cref{lem:root-moment}; the B\'ezoutian identity then gives
\begin{equation}\label{eq:Hermite-Bez}
B_{p,e}(x)[u]=C(x)H_{p,e}(x)[u]C(x)^\intercal.
\end{equation}
The identity extends to repeated roots by polynomial continuity in the coefficients of $a,b$. For the normal cubic on $t=0$, the matrix and its polynomial inverse are
\begin{equation}\label{eq:cubic-HB-matrix}
C(x)=\begin{pmatrix}-3r^2&0&1\\0&1&0\\1&0&0\end{pmatrix},
\qquad
C(x)^{-1}=\begin{pmatrix}0&0&1\\0&1&0\\1&0&3r^2\end{pmatrix}.
\end{equation}
For $p(e)\ne1$, normalize $p$ by $p(e)$; multiplying the resulting $C$ by $p(e)$ gives the congruence for $p$. In all cases both congruence matrices have polynomial entries, so matrix SOS is equivalent for the two parameterizations.

Our calculation was performed on the hyperplane transverse to $e_0$. Saunderson's translation congruence \cite[Proposition~3.11]{Saunderson2019} gives a unimodular polynomial matrix $K(t_0)$ satisfying
\begin{equation}\label{eq:translation-cong}
B_{p,e}(x+t_0e)[u]=K(t_0)B_{p,e}(x)[u]K(t_0)^\intercal .
\end{equation}
Thus being a matrix SOS on a transverse hyperplane implies being a matrix SOS in all variables.

The remaining coordinate changes and restrictions are covered by the following pullback identity.

\begin{lemma}\label{lem:linear-pullback}
Let $p$ be homogeneous and hyperbolic with respect to $e\in\R^n$. Let $A:\R^N\to\R^n$ be linear and let $\tilde e\in\R^N$ satisfy $A\tilde e=e$. Put $\tilde p=p\circ A$. Then
\[
\Lam(\tilde p,\tilde e)=A^{-1}(\Lam(p,e))
\]
and
\begin{equation}\label{eq:pullback-bez}
B_{\tilde p,\tilde e,\tilde u}(\tilde x)
=B_{p,e,A\tilde u}(A\tilde x).
\end{equation}
Consequently, SOS-hyperbolicity is preserved under such pullbacks. Positive rescaling of $p$ also preserves SOS-hyperbolicity.
\end{lemma}

\begin{proof}
For every $\tilde x$,
\[
\tilde p(\tilde x+t\tilde e)=p(A\tilde x+te),
\]
so the univariate roots agree and the cone identity follows. Likewise
\[
D_{\tilde u}\tilde p(\tilde x+t\tilde e)=D_{A\tilde u}p(A\tilde x+te),
\]
which proves \eqref{eq:pullback-bez}. Substitution preserves matrix SOS decompositions. If $p$ is replaced by $cp$ with $c>0$, its B\'ezoutian is multiplied by $c^2$.
\end{proof}

\begin{lemma}\label{lem:lower-m}
For every $1\le m\le4$ and $q\in\calH_m$, the polynomial
\[
p_q(t,x)=t^3-3t\|x\|^2+2q(x)
\]
is SOS-hyperbolic with respect to $(1,0)$.
\end{lemma}

\begin{proof}
The case $m=4$ follows from \cref{cor:H-matrix-sos,eq:Hermite-Bez,eq:translation-cong}. If $m<4$, define on $\R^4$ the cubic $\widetilde q(x,y)=q(x)$. Then
\[
|\widetilde q(x,y)|=|q(x)|\le\|x\|^3\le\|(x,y)\|^3,
\]
so $\widetilde q\in\calH_4$. The corresponding four-dimensional normal form
\[
\widetilde p(t,x,y)=t^3-3t\|(x,y)\|^2+2q(x)
\]
is SOS-hyperbolic. Restriction to the linear subspace $y=0$ gives $p_q$, and \cref{lem:linear-pullback} applies.
\end{proof}

\begin{theorem}\label{thm:five-variable-final}
Let $p$ be a homogeneous cubic in five real variables. For every hyperbolicity direction $e$ of $p$ and every $v\in\Lam(p,e)$, $B_{p,e,v}$ is a matrix SOS. Thus $p$ is SOS-hyperbolic in Saunderson's fixed-direction sense with respect to each such $e$.
\end{theorem}

\begin{proof}
Let $p$ be hyperbolic with respect to $e$. Normalize by the positive scalar $p(e)^{-1}$ and apply Netzer's normal-form reduction. If the transverse quotient has dimension zero, the reduced polynomial is $p(t)=t^3$; this is precisely the case where the original $p$ is a scalar multiple of the cube of a linear form. Retaining the transverse variables as constant directions, write $p(t,x)=t^3$ and $e_0=(1,0)$. A cone direction $u=(a,b)$ satisfies $a\ge0$, and direct calculation gives
\begin{equation}\label{eq:pure-cube-bez}
B_{t^3,e_0,u}(t,x)=3a\,c(t)c(t)^\intercal,
\qquad c(t)=(t^2,2t,1)^\intercal.
\end{equation}
Indeed, the defining divided difference is $3a(t+\tau)^2(t+s)^2$. Thus the B\'ezoutian is a rank-one matrix SOS when $a>0$, and is zero when $a=0$.

If the transverse quotient has positive dimension, the reduction gives $p_q$ with $q\in\calH_m$ and $1\le m\le4$, and \cref{lem:lower-m} applies. In both cases the inverse coordinate change and the discarded constant directions are linear pullbacks, handled by \cref{lem:linear-pullback}; positive rescaling restores the original normalization.
\end{proof}

\begin{remark}\label{rem:cubic-scope}
The argument uses no smoothness, irreducibility, or strict-hyperbolicity hypothesis. It is specific to transverse dimension at most four. Netzer proves that the set of Clifford-realizable cubics is a proper subset of the norm ball in transverse dimension at least five \cite[Theorem~2.5]{Netzer2026}. Thus the canonical-compression mechanism used here cannot extend unchanged to all cubic normal forms in six or more variables.
\end{remark}

\section{Two-direction weak--strong separation}\label{sec:separation}

We now construct the weak--strong separation. Both directions range over the closed hyperbolicity cone, as in the two-direction definitions of \cref{sec:prelim}. We first isolate a general localization principle, then verify its hypotheses for a rational cubic and transfer the obstruction to the final product. An interior perturbation gives the fixed-direction corollary.

\subsection{A general localization principle}

Localization uses two facts about products: a Wronskian identity for the scalar certificates and a cancellation property for square factors in the matrix transfer.

\begin{lemma}\label{lem:product-W}
For homogeneous polynomials $f,g$,
\begin{equation}\label{eq:Wfg}
W_{u,v}(fg)=f^2W_{u,v}(g)+g^2W_{u,v}(f).
\end{equation}
If $L=\prod_{k=0}^{N}\ell_k$ is a product of linear forms, then
\begin{equation}\label{eq:WL}
W_{u,v}(L)=\sum_{k=0}^{N}
\left(\frac{L}{\ell_k}\right)^2\ell_k(u)\ell_k(v).
\end{equation}
\end{lemma}

\begin{proof}
Expand the two directional derivatives in the definition of $W_{u,v}$. The mixed terms cancel, giving \eqref{eq:Wfg}. Apply \eqref{eq:Wfg} inductively to the product of linear factors and use $D_uD_v\ell_k=0$.
\end{proof}

\begin{lemma}\label{lem:square-cancel}
Let $B(z)$ be a real symmetric polynomial matrix and let $A(z)=z_1\cdots z_N$. If $A(z)^2B(z)$ is a matrix SOS, then $B(z)$ is a matrix SOS.
\end{lemma}

\begin{proof}
Write $A^2B=H^\intercal H$. Setting $z_1=0$ gives $H(0)^\intercal H(0)=0$. Each diagonal entry is an SOS equal to the zero polynomial, hence every entry of $H$ vanishes on the hyperplane $z_1=0$. Thus $z_1$ divides every entry of $H$. Write $H=z_1H_1$ and cancel $z_1^2$. Repeat for $z_2,\ldots,z_N$.
\end{proof}

\begin{proposition}\label{prop:localization}
Let $h\in\R[t,z_1,\ldots,z_N]_d$, with $N\ge1$, be homogeneous and hyperbolic with respect to $e=(1,0)$. Suppose that
\begin{enumerate}[label=(\roman*)]
\item $W_{e,e}(h)$ lies in the interior of the SOS cone of forms of degree $2d-2$ in $(t,z)$;
\item $B_{h,e,e}(0,z)$ is not a matrix SOS.
\end{enumerate}
For every sufficiently small $\eps>0$, set
\[
r_0=e,\qquad r_k=e+\eps e_k\ (1\le k\le N),\qquad
K_\eps=\operatorname{cone}(r_0,\ldots,r_N),
\]
where $e_k$ denotes the $k$th transverse coordinate vector. Define
\[
\ell_0=\eps t-\sum_{k=1}^N z_k,\qquad \ell_k=z_k\ (1\le k\le N),
\qquad L=\prod_{k=0}^N\ell_k,\qquad P=Lh.
\]
Then $P$ has degree $d+N+1$, is hyperbolic with respect to
$\bar e=\sum_{k=0}^N r_k$, and has cone $\Lam(P,\bar e)=K_\eps$.
It is two-direction weakly SOS-hyperbolic but not two-direction SOS-hyperbolic. If $h$ has rational coefficients, $\eps$ can be chosen rational, giving rational $P$ and $\bar e$.
\end{proposition}

The transverse formulation of hypothesis~(ii) is equivalent to failure of matrix SOS for $B_{h,e,e}$ in all variables, since \eqref{eq:translation-cong} extends a transverse factorization polynomially along $e$.

\begin{proof}
For sufficiently small $\eps$, all $r_k$ lie in $\operatorname{int}\Lam(h,e)$. By hypothesis (i), continuity and the finite number of generator pairs imply that every $W_{r_i,r_j}(h)$ is an SOS. Bilinearity then gives the same conclusion for all pairs in $K_\eps$.

To preserve the matrix obstruction, note that $B_{h,e,r_1}(0,z)\to B_{h,e,e}(0,z)$ in coefficients. Matrix SOS cones with a fixed degree bound are closed: a symmetric polynomial matrix $M(z)$ is a matrix SOS exactly when $y^\intercal M(z)y$ is scalar SOS in independent variables $(z,y)$, and the bounded-degree scalar SOS cone is closed. Indeed, in an SOS of this scalarization, homogeneity of degree two in $y$ forces each squared polynomial to be linear in $y$. Hypothesis (ii) therefore ensures that $B_{h,e,r_1}(0,z)$ is not a matrix SOS when $\eps$ is sufficiently small.

The half-space description
\[
K_\eps=\{(t,z):z_k\ge0,\ \eps t\ge\textstyle\sum_k z_k\}
\]
follows by solving for the coefficients of its rays. Hence $K_\eps\subseteq\Lam(h,e)$ and every $\ell_k$ is positive at $\bar e$. The product cone rule gives $\Lam(P,\bar e)=K_\eps$. For $u,v\in K_\eps$, \cref{lem:product-W} gives
\[
W_{u,v}(P)=L^2W_{u,v}(h)+h^2\sum_{k=0}^N
\left(\frac{L}{\ell_k}\right)^2\ell_k(u)\ell_k(v),
\]
an SOS. This proves two-direction weak SOS-hyperbolicity.

To transfer the obstruction, take $v=r_1$, put $A=\prod_{k=1}^Nz_k$ and $\alpha=-\sum_{k=1}^Nz_k$, and write $\ell(\tau)=\alpha+\eps\tau$. Since $\ell_0(v)=0$ and $A$ is constant along $e$,
\[
P(\tau,z)=A\ell(\tau)h(\tau,z),\qquad
D_vP=(D_vA)\ell h+A\ell D_vh.
\]
The term proportional to $\ell h$ contributes zero to the B\'ezoutian. Its nonzero $(d+1)\times(d+1)$ coefficient block is therefore
\begin{equation}\label{eq:general-block-transfer}
A^2T B_{h,e,v}(0,z)T^\intercal,
\end{equation}
where $T$ represents multiplication by $\ell$: $T_{j,j}=\alpha$ and $T_{j+1,j}=\eps$ for $0\le j<d$, with other entries zero. It has a polynomial left inverse $S$, given with indices $0\le i<d$, $0\le k\le d$ by
\begin{equation}\label{eq:general-left-inverse}
S_{ik}=\begin{cases}
(-\alpha)^{k-i-1}\eps^{-(k-i)},&k\ge i+1,\\
0,&k\le i.
\end{cases}
\qquad ST=I_d.
\end{equation}
Here $\eps$ is a fixed nonzero scalar. If $B_{P,e,v}$ were matrix SOS, its principal block \eqref{eq:general-block-transfer}, followed by congruence by $S$, would make $A^2B_{h,e,v}(0,z)$ matrix SOS. By \cref{lem:square-cancel}, $B_{h,e,v}(0,z)$ would then be matrix SOS, a contradiction. The directions $e=r_0$ and $v=r_1$ belong to $\partial K_\eps$. In particular $P(e)=0$, so $e$ is an admissible boundary direction in the two-direction condition, although it is not a hyperbolicity direction for $P$. Finally choose rational $\eps$ in the required open interval.
\end{proof}

This separates the seed conditions from the subsequent cone restriction and matrix transfer, which do not depend on the graph construction or on the cubic degree. We next supply a rational seed and a certified value of $\eps$.

\subsection{The icosahedral seed cubic}

For the explicit seed, let $G=(V,E)$ be the icosahedral graph on $V=\{1,\ldots,12\}$ with edge set
\begin{align}
E=\{&(1,2),(1,6),(1,8),(1,9),(1,12),(2,3),(2,6),(2,7),(2,9),\nonumber\\
&(3,4),(3,7),(3,9),(3,10),(4,5),(4,7),(4,10),(4,11),\nonumber\\
&(5,6),(5,7),(5,11),(5,12),(6,7),(6,12),(8,9),(8,10),\nonumber\\
&(8,11),(8,12),(9,10),(10,11),(11,12)\}.\label{eq:icosa-edges}
\end{align}
Thus $|V|=12$ and $|E|=30$. Introduce
\[
z=(x_1,\ldots,x_{12},(y_{ij})_{\{i,j\}\in E})\in\R^{42}
\]
and define
\begin{equation}\label{eq:Qg}
Q(z)\coloneqq \sum_{i=1}^{12}x_i^2+\sum_{\{i,j\}\in E}y_{ij}^2,
\qquad
g(z)\coloneqq \sum_{\{i,j\}\in E}x_ix_jy_{ij}.
\end{equation}
Set
\begin{equation}\label{eq:lambda}
\lambda\coloneqq \frac{199997}{200000}
\end{equation}
and
\begin{equation}\label{eq:seed-h}
h(t,z)\coloneqq t^3-3tQ(z)+9\lambda g(z),
\qquad e=(1,0,\ldots,0).
\end{equation}
The clique number of the icosahedral graph is $3$. The graph construction used by Saunderson gives
\begin{equation}\label{eq:gmax}
\max_{Q(z)=1}|g(z)|=\frac29;
\end{equation}
see \cite[Theorem~5.2]{Saunderson2019}. Hence the cubic $q_3(z)=\frac{9\lambda}{2}g(z)$ satisfies $\|q_3\|_\infty=\lambda<1$, so \eqref{eq:seed-h} is hyperbolic with respect to $e$.

The central Wronskian, obtained by direct differentiation, is
\begin{equation}\label{eq:Fseed}
W_{e,e}(h)=3F,
\qquad
F(t,z)\coloneqq t^4+3Q(z)^2-18\lambda t g(z).
\end{equation}

\subsection{An exact interior-SOS certificate}

We verify the first localization hypothesis with a positive-definite rational Gram matrix. The symmetries of the seed reduce its specification to finitely many orbit values.

Let $m_2$ be the column vector of all degree-two monomials in the $43$ variables $(t,z)$, in a fixed order. It has $\binom{44}{2}=946$ entries. Let
\[
\Gamma=(\mathbb Z/2\mathbb Z)^{13}
\]
act by independent signs $\eta_0,\eta_1,\ldots,\eta_{12}$ according to
\begin{equation}\label{eq:sign-action}
t\mapsto\eta_0t,
\qquad
x_i\mapsto\eta_ix_i,
\qquad
y_{ij}\mapsto\eta_0\eta_i\eta_j y_{ij}.
\end{equation}
The quartic $F$ is invariant under this action and under $\operatorname{Aut}(G)$, whose order is $120$.

Associate to each variable a character in $\mathbb F_2^{13}$ by
\[
\chi(t)=e_0,
\qquad
\chi(x_i)=e_i,
\qquad
\chi(y_{ij})=e_0+e_i+e_j,
\]
and extend additively to monomials. Define a symmetric rational matrix $G_F\in\Q^{946\times946}$ by the following rules:
\begin{enumerate}[label=(\roman*)]
\item $(G_F)_{\mu,\nu}=0$ if $\chi(\mu)\ne\chi(\nu)$;
\item for equal characters, the entry depends only on the $\operatorname{Aut}(G)$-orbit of the unordered pair $\{\mu,\nu\}$;
\item the $51$ orbit values are those in \cref{tab:orbit-data}.
\end{enumerate}
There are $2449$ unordered equal-character monomial pairs and $51$ automorphism orbits, so these rules determine the matrix.

\begin{proposition}\label{prop:GF}
The rational matrix $G_F$ satisfies
\begin{equation}\label{eq:GF}
F=m_2^\intercal G_Fm_2,
\qquad
G_F\succ0.
\end{equation}
Consequently $F$ lies in the interior of the quartic SOS cone, and $W_{e,e}(h)=3F$ is an interior SOS quartic.
\end{proposition}

\begin{proof}
Coefficient comparison using \cref{tab:orbit-data} establishes the identity in \eqref{eq:GF}. The sign-character decomposition block-diagonalizes $G_F$ into $564$ blocks. Under $\operatorname{Aut}(G)$ these blocks fall into the $15$ types listed in \cref{app:block-certificates}. The block dimensions account for every quadratic monomial:
\[
1\cdot43+12\cdot6+30\cdot4+30\cdot3+6+60+60\cdot2+15+60+60+60+60+20\cdot3+30\cdot2+60=946.
\]
For each block type, \cref{app:block-certificates} gives a positive integer scaling and the exact characteristic polynomial of the scaled symmetric block. Every displayed factor has alternating signs and nonzero constant term. Substituting $X=-u$ with $u>0$ shows that no factor has a nonpositive real root. Each block is real symmetric, so its eigenvalues are real and therefore strictly positive. Therefore $G_F\succ0$.
The rational coefficient and positive-definiteness checks are
reproducible as described in \cref{app:compute-gram}.
\end{proof}

\subsection{An exact non-SOS compression of the B\'ezoutian}

For the second localization hypothesis, we compress the central B\'ezoutian to a quartic and derive a contradiction from an exact trace bound.

At $t=0$ the cubic B\'ezoutian in the pair $(e,e)$ is
\begin{equation}\label{eq:Bseed}
B_{h,e,e}(0,z)=
\begin{pmatrix}
9Q^2&-27\lambda g&-3Q\\
-27\lambda g&6Q&0\\
-3Q&0&3
\end{pmatrix}.
\end{equation}
For an auxiliary scalar variable $s$, let
\begin{equation}\label{eq:Rseed}
\xi(s,z)\coloneqq (1,s,Q(z))^\intercal ,
\qquad
R(s,z)\coloneqq Q(z)^2+s^2Q(z)-9\lambda s g(z).
\end{equation}
Then
\begin{equation}\label{eq:Rcompression}
\xi(s,z)^\intercal B_{h,e,e}(0,z)\xi(s,z)=6R(s,z).
\end{equation}

The trace certificate uses the $20$ triangles of the icosahedral graph:
\begin{align}
\mathcal T=\{&(1,2,6),(1,2,9),(1,6,12),(1,8,9),(1,8,12),
(2,3,7),(2,3,9),(2,6,7),\nonumber\\
&(3,4,7),(3,4,10),(3,9,10),(4,5,7),(4,5,11),(4,10,11),
(5,6,7),(5,6,12),\nonumber\\
&(5,11,12),(8,9,10),(8,10,11),(8,11,12)\}.\label{eq:triangles}
\end{align}
Let
\[
m_0=(x_1^2,\ldots,x_{12}^2,(y_{ij}^2)_{\{i,j\}\in E})^\intercal \in\R^{42}.
\]
For each $C\in\mathcal T$, define $c_C\in\Q^{42}$ by
\begin{equation}\label{eq:cC}
(c_C)_{x_i^2}=\begin{cases}2/9,&i\in C,\\0,&i\notin C,\end{cases}
\qquad
(c_C)_{y_{ij}^2}=\begin{cases}1/9,&\{i,j\}\subset C,\\0,&\text{otherwise}.
\end{cases}
\end{equation}
Set
\begin{equation}\label{eq:SB}
S\coloneqq \sum_{C\in\mathcal T}c_Cc_C^\intercal ,
\qquad
B\coloneqq \diag\left(-I_{12},\frac{19}{50}I_{30}\right).
\end{equation}

\begin{lemma}\label{lem:trace-cert}
One has \[
B+802S\succ0.
\]
More precisely, the characteristic polynomial of $4050(B+802S)$ equals
\begin{align}
&(X-121839)^4(X-41639)^4(X-1539)^{10}\nonumber\\
&\quad\cdot(X^2-2644089X+2722171050)\nonumber\\
&\quad\cdot(X^2-639089X+84713850)^5\nonumber\\
&\quad\cdot\bigl(X^4-1839578X^3+196249414621X^2\nonumber\\
&\hspace{34mm}-167130585918900X+10172190802500\bigr)^3.\label{eq:trace-charpoly}
\end{align}
\end{lemma}

\begin{proof}
The characteristic polynomial is obtained by exact arithmetic from the rational vectors \eqref{eq:cC}. Each factor in \eqref{eq:trace-charpoly} has alternating signs and nonzero constant term. Since the matrix is symmetric, all roots are real; evaluation at $X=-u$, $u>0$, excludes nonpositive roots. Thus all eigenvalues are positive.
The exact matrix and characteristic-polynomial calculations are
reproduced in \cref{app:compute-trace}.
\end{proof}

\begin{proposition}\label{prop:R-not-sos}
The quartic $R$ in \eqref{eq:Rseed} is not an SOS. Consequently $B_{h,e,e}$ is not a matrix SOS.
\end{proposition}

\begin{proof}
Suppose $R$ were an SOS. Averaging its Gram matrix over the sign group \eqref{eq:sign-action}, with the auxiliary variable $s$ transforming as $t$, preserves positive semidefiniteness and makes the Gram matrix block diagonal by character. Since $R$ has no $s^4$ term, the diagonal entry corresponding to the monomial $s^2$ is zero. Positive semidefiniteness forces the entire row and column of $s^2$ to vanish.

The trivial-character block is therefore a PSD matrix $X$ on the basis $m_0$. The coefficients of $x_i^4$ and $y_{ij}^4$ in $R$ are all $1$, so
\begin{equation}\label{eq:diagX}
\diag X=\mathbf 1.
\end{equation}
The remaining character blocks contribute nonnegative values at every real point.

For a triangle $C$, set $s=1$, put $x_i=\sqrt{2/9}$ on the three vertices of $C$ and $0$ elsewhere, and put $y_{ij}=1/3$ on the three edges of $C$ and $0$ elsewhere. Then $Q=1$, $g=2/9$, and $m_0=c_C$. Therefore
\[
c_C^\intercal Xc_C\le R=2(1-\lambda).
\]
Summing over the $20$ triangles gives
\begin{equation}\label{eq:upperSX}
\tr(SX)\le40(1-\lambda)=\frac{3}{5000}.
\end{equation}
The opposite bound follows from \cref{lem:trace-cert}, $X\succeq0$, and \eqref{eq:diagX}:
\begin{align*}
0&\le\tr((B+802S)X)\\
&=\tr(BX)+802\tr(SX)\\
&=\left(-12+30\cdot\frac{19}{50}\right)+802\tr(SX)\\
&=-\frac35+802\tr(SX).
\end{align*}
Hence
\begin{equation}\label{eq:lowerSX}
\tr(SX)\ge\frac{3}{4010},
\end{equation}
contradicting \eqref{eq:upperSX}. Thus $R$ is not an SOS. If $B_{h,e,e}(0,z)$ were a matrix SOS, the polynomial compression \eqref{eq:Rcompression} would make $R$ an SOS. Hence the B\'ezoutian is not a matrix SOS.
\end{proof}

\subsection{A small rational cone}

Having verified the two seed hypotheses, we choose a cone on which all generator-pair Wronskians remain SOS. The sign symmetry will preserve the matrix obstruction under the required perturbation.

For $\eps>0$, define
\begin{equation}\label{eq:rays}
r_0\coloneqq e,
\qquad
r_k\coloneqq e+\eps e_k\quad(1\le k\le42),
\qquad
K_\eps\coloneqq \operatorname{cone}(r_0,\ldots,r_{42}),
\end{equation}
where $e_k$ is the $k$th coordinate vector in the $z$-space.

\begin{lemma}\label{lem:small-eps}
The following statements hold:
\begin{enumerate}[label=(\roman*)]
\item for all sufficiently small $\eps>0$, $r_0,\ldots,r_{42}\in\operatorname{int}\Lam(h,e)$;
\item for all sufficiently small $\eps>0$, every $W_{r_i,r_j}(h)$ is an SOS;
\item $B_{h,e,e+\eps e_1}(0,z)$ is not a matrix SOS for every real $\eps$, without a smallness restriction.
\end{enumerate}
In particular, there exists $\eps_0>0$ such that all three conclusions hold for $0<\eps<\eps_0$, and $\eps$ may be chosen rational.
\end{lemma}

\begin{proof}
Part (i) follows from openness of the cone interior and $r_k\to e$. For part (ii), \cref{prop:GF} places $W_{e,e}(h)$ in the interior of the finite-dimensional quartic SOS cone. The map $(u,v)\mapsto W_{u,v}(h)$ is bilinear and continuous in coefficients. There are finitely many pairs $(i,j)$, and $W_{r_i,r_j}(h)\to W_{e,e}(h)$ as $\eps\downarrow0$.

For part (iii), define
\[
R_\eps(s,z)\coloneqq \frac16\xi(s,z)^\intercal B_{h,e,e+\eps e_1}(0,z)\xi(s,z).
\]
Let $\sigma$ fix $s,t$ and all variables except $x_1$ and the incident edge variables $y_{1j}$, whose signs it reverses. Each summand of $g$ is unchanged, so
\[
h\circ\sigma=h,\qquad Q\circ\sigma=Q,\qquad \xi\circ\sigma=\xi.
\]
Moreover $\sigma e=e$ and $\sigma(e+\eps e_1)=e-\eps e_1$. Equivariance of directional differentiation and of the B\'ezoutian therefore gives
\[
R_\eps\circ\sigma=R_{-\eps}.
\]
Linearity in the second direction implies
\[
R=\frac12(R_\eps+R_{-\eps})
=\frac12(R_\eps+R_\eps\circ\sigma).
\]
If $R_\eps$ were SOS, substitution by $\sigma$ and averaging would make $R$ SOS, contradicting \cref{prop:R-not-sos}. Thus $R_\eps$ is not an SOS for any real $\eps$. Its defining compression excludes a matrix SOS decomposition of $B_{h,e,e+\eps e_1}(0,z)$ for every such $\eps$.

Intersecting the neighborhoods from (i) and (ii) completes the proof, since (iii) imposes no restriction on $\eps$. Density of $\Q$ in $\R$ then gives a rational choice of $\eps$.
\end{proof}

An explicit choice is $\eps=10^{-12}$, certified in \cref{app:compute-epsilon} by checking the $\binom{44}{2}=946$ unordered pairs of the $43$ cone generators, including repetitions. These are ray pairs, whereas the $2449$ equal-character pairs used to reconstruct $G_F$ are monomial pairs; the monomial vector has $946$ entries.

Fix henceforth $\eps=10^{-12}$, which satisfies \cref{lem:small-eps}. By bilinearity, if
\[
u=\sum_i a_ir_i,
\qquad
v=\sum_j b_jr_j,
\qquad a_i,b_j\ge0,
\]
then
\begin{equation}\label{eq:Wcone}
W_{u,v}(h)=\sum_{i,j}a_ib_jW_{r_i,r_j}(h),
\end{equation}
which is an SOS. Thus every Wronskian of $h$ with directions in $K_\eps$ is an SOS.

\subsection{Cutting the hyperbolicity cone by linear factors}

We now apply the product construction of \cref{prop:localization}. The linear factors impose the chosen polyhedral cone while preserving its scalar SOS certificates. Define
\begin{equation}\label{eq:linear-factors}
\ell_0(t,z)\coloneqq \eps t-\sum_{k=1}^{42}z_k,
\qquad
\ell_k(t,z)\coloneqq z_k\quad(1\le k\le42),
\end{equation}
and let
\begin{equation}\label{eq:final-P}
L(t,z)\coloneqq \prod_{k=0}^{42}\ell_k(t,z),
\qquad
P(t,z)\coloneqq L(t,z)h(t,z).
\end{equation}
Then $P$ is homogeneous in $43$ variables and has degree $46$.

The cone generated by \eqref{eq:rays} has the half-space description
\begin{equation}\label{eq:Khalf}
K_\eps=\left\{(t,z):z_k\ge0\ (1\le k\le42),\quad
\eps t\ge\sum_{k=1}^{42}z_k\right\}.
\end{equation}
Indeed, if $(t,z)=a_0r_0+\sum_{k=1}^{42}a_kr_k$ with $a_k\ge0$, then $z_k=\eps a_k$ and $\eps t=\eps a_0+\sum z_k$, and the converse follows by solving for the coefficients.

For the product, choose the hyperbolicity direction
\begin{equation}\label{eq:ebar}
\bar e\coloneqq \sum_{i=0}^{42}r_i=(43,\eps,\ldots,\eps).
\end{equation}
Then $\ell_k(\bar e)=\eps>0$ for all $k$. Since all $r_i$ lie in the interior of $\Lam(h,e)$, the cone $K_\eps$ is contained in $\Lam(h,e)$ and $\bar e$ lies in its interior. By the change-of-direction property, $h$ is hyperbolic with respect to $\bar e$ with the same cone. Each $\ell_k$ is hyperbolic with respect to $\bar e$ and has cone $\{\ell_k\ge0\}$. Therefore
\begin{equation}\label{eq:final-cone}
\Lam(P,\bar e)=\Lam(h,e)\cap\bigcap_{k=0}^{42}\{\ell_k\ge0\}=K_\eps.
\end{equation}

The scalar part of the localization argument now gives the following certificate; the matrix transfer is treated in the next subsection.

\begin{proposition}\label{prop:weak-final}
The polynomial $P$ in \eqref{eq:final-P} is two-direction weakly SOS-hyperbolic: $W_{u,v}(P)$ is an SOS for every $u,v\in\Lam(P,\bar e)$.
\end{proposition}

\begin{proof}
Let $u,v\in\Lam(P,\bar e)=K_\eps$. By \eqref{eq:Khalf}, $\ell_k(u),\ell_k(v)\ge0$. By \eqref{eq:Wcone}, $W_{u,v}(h)$ is an SOS. Therefore \cref{lem:product-W} gives
\begin{equation}\label{eq:WP-sos}
W_{u,v}(P)
=L^2W_{u,v}(h)
+h^2\sum_{k=0}^{42}\left(\frac{L}{\ell_k}\right)^2\ell_k(u)\ell_k(v),
\end{equation}
and every term on the right is an SOS.
\end{proof}

\subsection{Persistence of the matrix SOS obstruction}

The seed conditions in \cref{prop:localization} were established in \cref{prop:GF,prop:R-not-sos}. We now specialize the matrix transfer to $d=3$. Take the boundary directions
\[
u=e=r_0,
\qquad
v=r_1=e+\eps e_1,
\]
both in $\partial K_\eps$. Since $P(e)=0$, the first direction $e$ is not a hyperbolicity direction for $P$; it is admissible because the two-direction definition allows every direction in the closed cone. The hyperbolicity direction of the product remains $\bar e$. Since $\ell_0(v)=0$, define
\[
A(z)\coloneqq \prod_{k=1}^{42}z_k,
\qquad
\alpha(z)\coloneqq \ell_0(0,z)=-\sum_{k=1}^{42}z_k.
\]
Along the $e$-line through $(0,z)$,
\[
P(\tau,z)=A(z)(\alpha+\eps\tau)h(\tau,z).
\]
Moreover,
\[
D_vP(\tau,z)
=(D_vA)(z)(\alpha+\eps\tau)h(\tau,z)
+A(z)(\alpha+\eps\tau)D_vh(\tau,z).
\]
Over the coefficient ring in $z$, the first summand is proportional to the univariate polynomial $P$, so it contributes zero to the B\'ezoutian. Multiplication by $\alpha+\eps\tau$ acts on degree at most two coefficient vectors through
\begin{equation}\label{eq:Talpha}
T(\alpha)=
\begin{pmatrix}
\alpha&0&0\\
\eps&\alpha&0\\
0&\eps&\alpha\\
0&0&\eps
\end{pmatrix}.
\end{equation}
Consequently, the nonzero $4\times4$ coefficient block of $B_{P,e,v}(0,z)$ is
\begin{equation}\label{eq:block-transfer}
A(z)^2T(\alpha)B_{h,e,v}(0,z)T(\alpha)^\intercal .
\end{equation}
All remaining rows and columns of the degree-$46$ B\'ezoutian vanish because $P((0,z)+\tau e)$ has degree only four in $\tau$.

The matrix $T(\alpha)$ has the polynomial left inverse
\begin{equation}\label{eq:Salpha}
S(\alpha)=
\begin{pmatrix}
0&\eps^{-1}&-\alpha\eps^{-2}&\alpha^2\eps^{-3}\\
0&0&\eps^{-1}&-\alpha\eps^{-2}\\
0&0&0&\eps^{-1}
\end{pmatrix},
\qquad
S(\alpha)T(\alpha)=I_3.
\end{equation}

\begin{theorem}\label{thm:separation-final}
The polynomial $P$ in \eqref{eq:final-P} has rational coefficients, degree $46$, and $43$ variables. It is two-direction weakly SOS-hyperbolic: $W_{u,v}(P)$ is an SOS for all $u,v\in\Lam(P,\bar e)$. It is not two-direction SOS-hyperbolic, since $B_{P,e,r_1}$ is not a matrix SOS, with $e,r_1\in\partial\Lam(P,\bar e)$. No failure of the fixed-direction matrix SOS condition at $\bar e$ is asserted.
\end{theorem}

\begin{proof}
Two-direction weak SOS-hyperbolicity is \cref{prop:weak-final}. Suppose, for contradiction, that $B_{P,e,v}$ were a matrix SOS. Then the principal block \eqref{eq:block-transfer} would be a matrix SOS. Polynomial congruence by the left inverse \eqref{eq:Salpha} would imply that
\[
A(z)^2B_{h,e,v}(0,z)
\]
is a matrix SOS. By \cref{lem:square-cancel}, $B_{h,e,v}(0,z)$ would be a matrix SOS, contradicting \cref{lem:small-eps}(iii).
\end{proof}

The matrix obstruction can also be moved into the cone interior. This converts the two-direction separation into a fixed-direction separation without changing the polynomial.

\begin{samepage}
\begin{corollary}\label{cor:fixed-direction-separation}
Let $P$, $K_\eps$, $e=r_0$, $r_1$, and $\bar e$ be as above. For every sufficiently small $\delta>0$, set
\begin{equation}\label{eq:interior-directions}
e_\delta=e+\delta\bar e,\qquad
v_\delta=r_1+\delta\bar e.
\end{equation}
Then $e_\delta,v_\delta\in\operatorname{int}K_\eps$, $P$ is hyperbolic with respect to $e_\delta$, and
\[
\Lam(P,e_\delta)=K_\eps,\qquad
W_{e_\delta,v}(P)\text{ is an SOS for every }v\in K_\eps,
\]
whereas $B_{P,e_\delta,v_\delta}$ is not a matrix SOS. Thus $P$ is weakly SOS-hyperbolic but not SOS-hyperbolic with respect to $e_\delta$ in Saunderson's fixed-direction sense. The parameter $\delta$ can be chosen rational, giving rational $e_\delta$ and $v_\delta$.
\end{corollary}
\end{samepage}

\begin{proof}
Each defining linear form of $K_\eps$ is nonnegative on $e,r_1$ and strictly positive on $\bar e$. Therefore $e_\delta$ and $v_\delta$ lie in $\operatorname{int}K_\eps$ for every $\delta>0$. The change-of-direction property gives hyperbolicity with respect to $e_\delta$ and the same cone. By \cref{prop:weak-final}, $W_{e_\delta,v}(P)$ is an SOS for every $v\in K_\eps$.

By \eqref{eq:parambez}, the coefficients of the parameterized B\'ezoutian depend polynomially on its two directions. Consequently,
\[
B_{P,e_\delta,v_\delta}\longrightarrow B_{P,e,r_1}
\qquad\text{in coefficients as }\delta\downarrow0.
\]
These matrices have size $46$ and entries of degree at most $90$ in the base variables. By the scalarization argument in \cref{prop:localization}, the matrix SOS cone with these bounds is closed, so its complement is open. Since $B_{P,e,r_1}$ is not a matrix SOS by \cref{thm:separation-final}, neither is $B_{P,e_\delta,v_\delta}$ for all sufficiently small $\delta>0$. This gives failure of the fixed-direction matrix SOS condition at $e_\delta$, witnessed by the interior cone direction $v_\delta$. Finally choose rational $\delta$ in the resulting open interval.
\end{proof}

\begin{remark}
\Cref{cor:fixed-direction-separation} disproves unconditional equivalence of the fixed-direction weak and strong notions in Saunderson's Question~4.6 \cite{Saunderson2019}.
\end{remark}

\section{Consequences and scope}\label{sec:discussion}

The results obtained in this paper distinguish the roles of auxiliary degree, cubic structure, and scalar versus matrix SOS certificates. We summarize their scope and the questions they leave open.

First, unrestricted hyperzout representations are existentially complete for nonnegative forms. The degree restriction in \eqref{eq:restricted-degree} is essential to the non-representability results: once the auxiliary degree is allowed to reach the boundary $dk=2r$, every nonnegative form can be represented. The construction does not determine the smallest auxiliary degree or the smallest ambient dimension for a given form, so these minimization questions are not addressed here.

Second, fixed-direction exactness for five-variable cubics holds for each hyperbolicity direction $e$: for every $v\in\Lam(p,e)$, $B_{p,e,v}$ is a matrix SOS. Yet, the present argument is restricted to transverse dimension at most four.

Third, \cref{prop:localization} separates scalar and matrix certificates: all $W_{u,v}(P)$ on the final cone are SOS, but $B_{P,e,r_1}$ is not a matrix SOS. The icosahedral seed satisfies the hypotheses, and its sign symmetry preserves the compressed obstruction for every real perturbation parameter. Although $e=r_0$ is a boundary direction with $P(e)=0$, \cref{cor:fixed-direction-separation} gives fixed-direction separation at a rational hyperbolicity direction $e_\delta$. This disproves unconditional equivalence in Saunderson's Question~4.6. The matrix SOS status at $\bar e$ remains undetermined, and separation under irreducibility remains open.

Two further directions follow. Quantitatively, one can seek sharp auxiliary-degree and dimension bounds for hyperzout representations of specified classes of nonnegative forms. Structurally, one can ask when pointwise B\'ezoutian positivity, scalar SOS conditions, or spectrahedrality imply matrix SOS decompositions.

\section*{Acknowledgements}
The authors acknowledge the use of GPT-6 Astra to assist with
brainstorming, mathematical development, and manuscript drafting.
The authors are solely responsible for the final content, analysis, and
conclusions.

\section*{Funding}
The authors were supported by grant 2023YFA1009401 from the National Key R\&D Program of China and grant 12171324 from the National Natural Science Foundation of China.

\bibliographystyle{amsplain}
\bibliography{hyperbolic_polynomial}

\newpage
\appendix
\section{Exact rational certificates for the separation theorem}\label{app:certificates}

This appendix specifies the rational Gram matrix used in \cref{prop:GF} through its orbit values and scaled block characteristic polynomials. The final subsection describes how to reconstruct and verify these data exactly.

\subsection{\texorpdfstring{Orbit data for $G_F$}{Orbit data for G F}}

In \cref{tab:orbit-data}, $y_{ij}$ denotes the edge variable corresponding to $\{i,j\}\in E$. Each row gives a representative unordered pair of degree-two monomials and the common rational entry assigned to its full $\operatorname{Aut}(G)$-orbit. Entries joining unequal sign characters are zero.

\begin{longtable}{r>{\raggedright\arraybackslash}p{0.46\textwidth}>{\raggedleft\arraybackslash}p{0.31\textwidth}}
\caption{Orbit data defining $G_F$.}\label{tab:orbit-data}\\
\toprule
$k$ & representative $\{\mu,\nu\}$ & $q_k$\\
\midrule
\endfirsthead
\toprule
$k$ & representative $\{\mu,\nu\}$ & $q_k$\\
\midrule
\endhead
0 & $\{x_1y_{34},x_1y_{34}\}$ & $277703341/50000000$\\
1 & $\{x_1y_{26},x_2y_{16}\}$ & $-4704751/1562500$\\
2 & $\{y_{12}y_{34},y_{12}y_{34}\}$ & $187424861/25000000$\\
3 & $\{x_1y_{12},x_3y_{23}\}$ & $127435361/20000000$\\
4 & $\{x_1^2,y_{23}^2\}$ & $172400307/200000000$\\
5 & $\{y_{12}y_{37},y_{12}y_{37}\}$ & $77665897/10000000$\\
6 & $\{y_{12}y_{39},y_{19}y_{23}\}$ & $-423175173/100000000$\\
7 & $\{y_{12}y_{16},y_{12}y_{16}\}$ & $13441099/2000000$\\
8 & $\{x_1x_3,y_{12}y_{23}\}$ & $-209630147/100000000$\\
9 & $\{y_{12}^2,y_{16}^2\}$ & $-1441099/4000000$\\
10 & $\{t^2,x_1^2\}$ & $-141593283/100000000$\\
11 & $\{y_{12}^2,y_{12}^2\}$ & $3$\\
12 & $\{ty_{12},y_{16}y_{26}\}$ & $0$\\
13 & $\{x_1y_{37},x_1y_{37}\}$ & $334885833/50000000$\\
14 & $\{y_{12}^2,y_{18}^2\}$ & $-68007329/200000000$\\
15 & $\{x_1^2,x_3^2\}$ & $48784799/25000000$\\
16 & $\{y_{12}y_{18},y_{12}y_{18}\}$ & $668007329/100000000$\\
17 & $\{x_1y_{12},x_1y_{12}\}$ & $318608967/50000000$\\
18 & $\{x_1^2,y_{37}^2\}$ & $-34885833/100000000$\\
19 & $\{y_{12}^2,y_{39}^2\}$ & $176718797/200000000$\\
20 & $\{x_1^2,y_{34}^2\}$ & $22296659/100000000$\\
21 & $\{x_1^2,x_2^2\}$ & $356336261/200000000$\\
22 & $\{y_{12}^2,y_{34}^2\}$ & $-37424861/50000000$\\
23 & $\{x_1y_{23},x_1y_{23}\}$ & $427599693/100000000$\\
24 & $\{x_1^2,y_{12}^2\}$ & $-18608967/100000000$\\
25 & $\{x_1y_{26},x_1y_{26}\}$ & $301107077/50000000$\\
26 & $\{y_{12}y_{39},y_{12}y_{39}\}$ & $423281203/100000000$\\
27 & $\{x_1x_3,x_1x_3\}$ & $26215201/12500000$\\
28 & $\{y_{12}y_{18},y_{29}y_{89}\}$ & $-248824291/100000000$\\
29 & $\{y_{12}y_{3,10},y_{12}y_{3,10}\}$ & $249137789/50000000$\\
30 & $\{ty_{12},ty_{12}\}$ & $3362051/10000000$\\
31 & $\{x_1x_2,y_{16}y_{26}\}$ & $-336020823/100000000$\\
32 & $\{x_1x_2,x_1x_2\}$ & $243663739/100000000$\\
33 & $\{x_1^2,y_{26}^2\}$ & $-1107077/100000000$\\
34 & $\{x_1y_{23},x_3y_{12}\}$ & $-213773329/50000000$\\
35 & $\{x_1y_{12},x_6y_{26}\}$ & $637124887/100000000$\\
36 & $\{t^2,y_{12}^2\}$ & $-3362051/20000000$\\
37 & $\{y_{12}^2,y_{3,10}^2\}$ & $50862211/100000000$\\
38 & $\{y_{12}^2,y_{45}^2\}$ & $109373/390625$\\
39 & $\{x_1^2,x_1^2\}$ & $3$\\
40 & $\{y_{12}y_{16},y_{27}y_{67}\}$ & $83999933/12500000$\\
41 & $\{y_{12}y_{45},y_{12}y_{45}\}$ & $2125004/390625$\\
42 & $\{y_{12}^2,y_{37}^2\}$ & $-17665897/20000000$\\
43 & $\{tx_1,x_2y_{12}\}$ & $-212389003/50000000$\\
44 & $\{t^2,t^2\}$ & $1$\\
45 & $\{ty_{12},x_1x_2\}$ & $-6303811/12500000$\\
46 & $\{tx_1,tx_1\}$ & $141593283/50000000$\\
47 & $\{y_{12}y_{4,11},y_{12}y_{4,11}\}$ & $112043547/25000000$\\
48 & $\{x_1^2,x_4^2\}$ & $121010621/50000000$\\
49 & $\{x_1x_4,x_1x_4\}$ & $28989379/25000000$\\
50 & $\{y_{12}^2,y_{4,11}^2\}$ & $37956453/50000000$\\
\bottomrule
\end{longtable}

\subsection{\texorpdfstring{Exact positivity certificate for $G_F$}{Exact positivity certificate for G F}}\label{app:block-certificates}

The character decomposition gives $564$ blocks in $15$ automorphism types. For each type below, ``mult.'' is the number of blocks of that type, ``size'' is the block size, and $d_j$ is a positive integer such that the displayed polynomial is the characteristic polynomial of $d_jG_j$.

\begin{enumerate}[label=\textbf{Type \arabic*.},leftmargin=4.5em]
\setcounter{enumi}{-1}
\item mult. $1$, size $43$, $d_0=200000000$:
\begin{align*}
&(X-854270790)^5(X-509460)^4(X-49542)^4\\
&\quad\cdot(X^3-5241027569X^2+305503575098265068X-18813190163221355560)\\
&\quad\cdot(X^3-2547715335X^2+632310126777944156X-4342386569824759941684)^5\\
&\quad\cdot(X^4-1048935266X^3+237398030060911143X^2\\
&\hspace{27mm}-59037556334482284251102X+112394644170406114092376124)^3.
\end{align*}
\item mult. $12$, size $6$, $d_1=100000000$:
\[
(X^2-3469007884X+27027333808)(X^2-134176X+1131451339)^2.
\]
\item mult. $30$, size $4$, $d_2=100000000$:
\[
(X-55486)(X^3-1621338663X^2+152513996875699174X-198879174346605309264).
\]
\item mult. $30$, size $3$, $d_3=100000000$:
\[
(X-916831620)(X^2-628904646X+22143713201886).
\]
\item mult. $6$, size $1$, $d_4=25000000$: $X-28989379$.
\item mult. $60$, size $1$, $d_5=50000000$: $X-277703341$.
\item mult. $60$, size $2$, $d_6=100000000$: $(X-855146351)(X-53035)$.
\item mult. $15$, size $1$, $d_7=25000000$: $X-112043547$.
\item mult. $60$, size $1$, $d_8=25000000$: $X-187424861$.
\item mult. $60$, size $1$, $d_9=50000000$: $X-334885833$.
\item mult. $60$, size $1$, $d_{10}=50000000$: $X-249137789$.
\item mult. $60$, size $1$, $d_{11}=390625$: $X-2125004$.
\item mult. $20$, size $3$, $d_{12}=50000000$: $(X-451659109)^2(X-3013)$.
\item mult. $30$, size $2$, $d_{13}=100000000$: $(X-846456376)(X-106030)$.
\item mult. $60$, size $1$, $d_{14}=10000000$: $X-77665897$.
\end{enumerate}

Every factor has alternating signs with nonzero constant term, and each block is real symmetric. This gives the exact positivity proof used in \cref{prop:GF}.

\subsection{Independent verification recipe}

The orbit and block data above give the following exact reconstruction procedure.
\begin{enumerate}
\item Build the graph from \eqref{eq:icosa-edges} and enumerate its automorphism group; its order is $120$.
\item Enumerate the $946$ degree-two monomials, compute their characters from \eqref{eq:sign-action}, and set all cross-character Gram entries to zero.
\item Enumerate automorphism orbits of the $2449$ unordered equal-character monomial pairs (including diagonal pairs), and assign the $51$ values in \cref{tab:orbit-data}. Expand $m_2^\intercal G_Fm_2$ and compare every coefficient with $F$ in \eqref{eq:Fseed}.
\item Split $G_F$ by sign character and compute exact characteristic polynomials of the blocks after the listed integer scalings. They must match \cref{app:block-certificates}.
\item Form the $20$ vectors $c_C$ from \eqref{eq:cC}, compute $S$, and compute the exact characteristic polynomial of $4050(B+802S)$. It must equal \eqref{eq:trace-charpoly}.
\end{enumerate}
% No semidefinite solver output or numerical tolerance is required in these checks.

% Appendix B: only arithmetic checks requiring large-scale computation.
% Existing Appendix A is retained verbatim.
\section{Reproducible large-scale exact computations for the separation theorem}
\label{app:computer}

The supplementary files\footnote{They can be downloaded at \url{https://wangjie212.github.io/jiewang/codes/hyperbolic_verification.rar}.} isolate the six computationally substantial claims
used in Section~5 and Appendix~A.1--A.2. The single input file
\texttt{paper\_thm3\_data.json} reproduces the graph edges in
\eqref{eq:icosa-edges}, the value of $\lambda$ in \eqref{eq:lambda},
and the $51$ orbit entries in \cref{tab:orbit-data}.
Each script can be run separately: it contains all of its own helper
routines, reads the same JSON input, and does not import or execute
any other verification script. With Python~3, NetworkX,
and SymPy installed, execute
\begin{verbatim}
python -m pip install -r requirements.txt
python verify_all_heavy.py
\end{verbatim}
Alternatively, run each \texttt{check\_NN\_*.py} separately or use the
fully combined \texttt{verify\_all\_in\_one.py}. All certificate
comparisons use exact integers or rational numbers; the characteristic
polynomials are computed symbolically over $\mathbb Z$.

\subsection{Reconstruction of the Gram certificate}
\label{app:compute-gram}

The file \texttt{check\_01\_gram\_identity.py} independently rebuilds the
$946\times946$ matrix $G_F$ from the $51$ orbit representatives in
\cref{tab:orbit-data}: $120$ graph automorphisms expand those representatives
into all $2449$ equal-character unordered pairs, arranged in $564$
character blocks. It then compares every quartic coefficient to establish
\begin{equation}\label{eq:compute-gram-identity}
F=m_2^\intercal  G_Fm_2.
\end{equation}
The file \texttt{check\_02\_gram\_positive.py} checks each of the
$564$ shifted blocks by exact integer Bareiss elimination and Sylvester's
criterion. It establishes the strict rational margin
\begin{equation}\label{eq:compute-gram-margin}
G_F-10^{-8}I_{946}\succ0.
\end{equation}
Both files concern Proposition~\ref{prop:GF}: agreement of the polynomial
coefficients and positivity are separate requirements for the claimed
SOS certificate.

The separate file \texttt{check\_03\_gram\_charpoly.py} compares all $15$
scaled block characteristic polynomials, including multiplicities, against
the large integer factorizations in \cref{app:block-certificates}. This
last check confirms the printed computational data; the direct
positive-definiteness check above already suffices to establish the
strict Gram certificate.

\subsection{The trace certificate}
\label{app:compute-trace}

The file \texttt{check\_04\_trace\_positive.py} constructs the $20$
triangle vectors $c_C$ from \eqref{eq:triangles}--\eqref{eq:cC}, forms
$S$ and $B$ as in \eqref{eq:SB}, and checks the strict positive
definiteness of the exact integer matrix $4050(B+802S)$. It therefore
establishes the assertion $B+802S\succ0$ in Lemma~\ref{lem:trace-cert}.
The separate file \texttt{check\_05\_trace\_charpoly.py} compares
the characteristic polynomial of that $42\times42$ matrix with the
complete factored expression in \eqref{eq:trace-charpoly}.
These two calculations are the only large matrix computations needed
for the proof of Proposition~\ref{prop:R-not-sos}; the subsequent trace
inequality and scalar SOS contradiction are proved in the main text.

\subsection{All generator-pair perturbations}
\label{app:compute-epsilon}

The remaining large computation is the verification, for one fixed
rational parameter, of all $\binom{44}{2}=946$ unordered pairs of
rays in \eqref{eq:rays}. The file
\texttt{check\_06\_946\_wronskians.py} forms the necessary directional
Wronskians of the cubic \eqref{eq:seed-h} and compares every coefficient
of $W_{r_i,r_j}(h)-W_{e,e}(h)$. Take
\begin{equation}\label{eq:compute-epsilon}
\eps=10^{-12}.
\end{equation}
The exact maximum coefficient $\ell^1$-norm is
\begin{equation}\label{eq:compute-epsilon-bound}
\max_{0\le i\le j\le42}
\left\|W_{r_i,r_j}(h)-W_{e,e}(h)\right\|_{\mathrm{coef}}
=\frac{573351414000184315458626973}
{40000000000000000000000000000000000}
<3\cdot10^{-8},
\end{equation}
attained at $(i,j)=(1,2)$. The coefficient norm is the sum of the
absolute values of the coefficients in the ordinary quartic monomial
basis. For every quartic $H$ one can split each quartic monomial into
a product of two quadratic monomials, thus obtaining a symmetric Gram
matrix $G_H$ such that
$H=m_2^\intercal  G_Hm_2$ and
$\|G_H\|_{\mathrm{op}}\le\|H\|_{\mathrm{coef}}$.
Hence \eqref{eq:compute-gram-margin},
$W_{e,e}(h)=3F$, and \eqref{eq:compute-epsilon-bound} imply that all
$946$ generator-pair Wronskians have positive-definite Gram matrices.
The bilinear extension to every $u,v\in K_\eps$ remains the elementary
convexity argument \eqref{eq:Wcone} in the main text.

The claim in Lemma~\ref{lem:small-eps}(iii) requires no massive
symbolic verification. Indeed, reverse the signs of $x_1$ and of all
$y_{1j}$ incident to vertex $1$. Each term of $h$ is unchanged, while
the direction $e+\eps e_1$ is mapped to $e-\eps e_1$. If
$B_{h,e,e+\eps e_1}(0,z)$ were a matrix SOS, substitution by this
involution would also give a matrix SOS representation for
$B_{h,e,e-\eps e_1}(0,z)$. Their average equals $B_{h,e,e}(0,z)$
by linearity in the second direction, contradicting
Proposition~\ref{prop:R-not-sos}. This short argument applies for every
$\eps>0$ and is not implemented as a separate Python check.
The inclusion $r_i\in\operatorname{int}\Lam(h,e)$ for this choice follows
from the simple formula
$h(\tau,\eps e_i)=\tau^3-3\eps^2\tau$ and
$\sqrt3\eps<1$.

\subsection{Scope of the calculation files}
\label{app:compute-scope}

We deliberately omit software checks of the elementary B\'ezoutian
formula \eqref{eq:Bseed}, the compression \eqref{eq:Rcompression},
the Wronskian product rule \eqref{eq:Wfg}, the factorization
\eqref{eq:block-transfer}, and the left-inverse identity
\eqref{eq:Salpha}. These are short symbolic calculations explicitly
explained in Section~5, not difficult data-dependent assertions.
The proof of Theorem~\ref{thm:separation-final} uses the large
computations above solely through the Gram certificate, the trace
certificate, and the all-pairs perturbation bound.

\ifprintcomputercode
\clearpage
\section{Complete executable large-scale verification sources}
\label{app:python-listing}
\noindent The following standalone combined file implements exactly
those six checks; each of the six per-claim scripts is supplied
separately as supplementary material.
\lstinputlisting{verify_all_in_one.py}
\fi

\end{document}